\documentclass[a4paper,11pt]{article}
\usepackage{authblk}
\usepackage[english]{babel}
\usepackage{amsfonts,amsthm,amssymb}
\usepackage{enumerate}
\usepackage{placeins}
\usepackage[hmargin={28mm,28mm},vmargin={30mm,35mm}]{geometry}
\usepackage{relsize}
\usepackage{xcolor}
\usepackage{tcolorbox}
\usepackage[ocgcolorlinks,allcolors={blue}]{hyperref}
\usepackage{nicefrac}
\usepackage{mathtools,mathrsfs}
\usepackage{bm}
\usepackage{bbm}
\usepackage{mathabx}
\usepackage{setspace}
\usepackage{subfig}
\usepackage{floatrow}
\usepackage{a4wide}
\usepackage{color}
\usepackage[applemac]{inputenc}
\usepackage[T1]{fontenc}
\usepackage{url}
\usepackage{upgreek}
\usepackage{paralist}
\usepackage{accents}
\allowdisplaybreaks
\usepackage{authblk}
\usepackage{parskip}
\usepackage{algorithmic}

\begingroup
    \makeatletter
    \@for\theoremstyle:=definition,remark,plain\do{%
        \expandafter\g@addto@macro\csname th@\theoremstyle\endcsname{%
            \addtolength\thm@preskip\parskip
            }%
        }
\endgroup

\usepackage{diagbox}
\usepackage{booktabs}

\makeatletter
\def\thm@space@setup{%
  \thm@preskip=\parskip \thm@postskip=0pt
}
\makeatother

\newcommand{\arxiv}[1]{\href{https://arxiv.org/abs/#1}{arXiv:#1}}

\newcommand{\cs}{c_{\mathrm{sch}}}

\newcommand{\bbbracketleft}{[\hspace{-2.2pt}[}
\newcommand{\bbbracketright}{]\hspace{-2.2pt}]}
\newcommand{\curlybraceleft}{\{\hspace{-3.2pt}\{}
\newcommand{\curlybraceright}{\}\hspace{-3.2pt}\}}
\newcommand{\jm}[1]{{\bbbracketleft{#1}\bbbracketright}}
\newcommand{\av}[1]{{\curlybraceleft{#1}\curlybraceright}}

\newcommand{\sF}[1]{\mathsf{#1}}

\newcommand{\NR}{{N\!R}}
\newcommand*{\dotp}{\raisebox{-0.2ex}
{\scalebox{1.3}{$\cdot$}}}
\renewcommand*{\ddot}[1]{%
  \accentset{\mbox{\footnotesize\bfseries .\hspace{-0.05ex}.}}{#1}}
 \renewcommand*{\dot}[1]{%
  \accentset{\mbox{\footnotesize\bfseries .}}{#1}}

\newcommand{\E}{{\cc E}}

\newcommand{\Del}{\triangle}

\newcommand{\DG}{{\mathrm{DG}}}

\newcommand{\wl}{\overline}
\newcommand{\la}{\left\langle}
\newcommand{\ra}{\right\rangle}

\newcommand{\bb}{\mathbb}
\newcommand{\wt}{\widetilde}
\let\temp\phi
\let\phi\varphi
\let\varphi\temp

\newcommand{\nf}{\nicefrac}
\newcommand{\cc}{\mathcal}
\newcommand{\de}{\partial}

\newcommand{\eps}{\varepsilon}
\newcommand{\wh}{\widehat}
\newcommand{\z}{\zeta}
\newcommand{\K}{{\cc K}}

\newcommand{\Oe}{\Omega_e}
\newcommand{\Oa}{\Omega_a}

\newcommand{\sym}    {\mbox{\rm sym}\,}
\newcommand{\grad}  {\bm\nabla}

\newcommand{\divb} {\mathbf{div}\,}		
\newcommand{\tr}    {\mbox{\rm tr}\,}

\renewcommand{\bf}{\mathbf}
\newcommand{\mr}{\mathrm}
\renewcommand{\d}{\mr{d}}

\newtheorem{theorem}{Theorem}
\numberwithin{theorem}{section}

\theoremstyle{remark}
\newtheorem{remark}[theorem]{Remark}
\theoremstyle{definition}


\newcounter{parentnumber}
\newcommand{\email}[1]{\href{mailto:#1}{#1}}

\begin{document}
\title{Simulation of 3D elasto-acoustic wave propagation based on {a} Discontinuous Galerkin Spectral Element {method}\footnote{This work has been supported by SIR Research Grant no.~RBSI14VTOS funded by MIUR -- Italian Ministry of Education, Universities, and Research, and by ``National Group of Scientific Computing'' (GNCS-INdAM).}}
\author[1]{{Paola F.~Antonietti}\footnote{\email{paola.antonietti@polimi.it}}}
\author[1]{{Francesco Bonaldi}\footnote{Corresponding author, 
\email{francesco.bonaldi@polimi.it}}}
\author[1]{{Ilario Mazzieri}\footnote{\email{ilario.mazzieri@polimi.it}}}
\affil[1]{MOX, Dipartimento di Matematica, Politecnico di Milano\\ Piazza Leonardo da Vinci 32, 20133 Milano, Italy}%
\date{}
\maketitle
\begin{abstract}
\noindent
In this paper we present a numerical discretization of the coupled elasto-acoustic wave propagation problem based on a Discontinuous Galerkin Spectral Element (DGSE) approach in a three-dimensional setting. The unknowns of the coupled problem are the displacement field and the velocity potential, in the elastic and the acoustic domains, respectively, thereby resulting in a symmetric formulation. After stating the main theoretical results, we assess the performance of the method by convergence tests carried out on both matching and non-matching grids, and we simulate {realistic} scenarios where elasto-acoustic coupling occurs. In particular, we consider the case of Scholte waves and the scattering of elastic waves by an underground acoustic cavity. Numerical simulations are carried out by means of the code \texttt{SPEED}, {available at \texttt{http://speed.mox.polimi.it}}.
\end{abstract}

\section*{Introduction}
The main goal of this work is to simulate three-dimensional scenarios of elasto-acoustic coupling via a Discontinuous Galerkin Spectral Element (DGSE) discretization. 
Coupled elasto-acoustic wave propagation arises in several scientific and engineering contexts. In a geophysical framework, a first example one can think of is given by seismic events occurring near coastal environments; another relevant situation where such a problem plays a major role is the detection of underground cavities \cite{mazzieri-cavity,perugia-cavity,perugia-cavity-2}. Elasto-acoustic coupling occurs in structural acoustics as well, when sensing or actuation devices are immersed in an acoustic fluid \cite{flemisch-wohlmuth}, and also in medical ultrasonics \cite{monkola-sanna,monkola-sanna-paper}.

Typically, an elasto-acoustic coupling arises in the following framework: a space region made up by two subregions, one occupied by a solid (elastic) medium, the other by a fluid (acoustic) one, with suitable transmission conditions imposed at the interface between the two. The aim of such  conditions is to account for the following physical properties:~\begin{inparaenum}[(i)] \item the normal component of the velocity field is continuous at the interface;
\item a pressure load is exerted by the fluid body on the solid one through the interface.
\end{inparaenum}
In a geophysics context, when a seismic event occurs near a coastal environment, both \emph{pressure} (P) and \emph{shear} (S) waves are generated. However, only P-waves (i.e., whose direction of propagation is aligned with the displacement of the medium) are able to travel through both solid and fluid media, unlike S-waves (i.e., whose direction of propagation is orthogonal to the displacement of the medium), which can travel only through solids. This explains the reason for considering the first interface condition. On the other hand, the second one accounts for the fact that an acoustic wave propagating in a fluid domain gives rise to an \emph{acoustic pressure} exerted on the solid via the interface.

Numerical simulation of elasto-acoustic coupling scenarios has been the subject of a very broad literature. We give below a brief (and by far non-exhaustive) overview of some of the research works carried out so far in this field. Bathe \emph{et al.}~\cite{bathe} and Berm\'udez \emph{et al.}~\cite{bermudez} considered a displacement-based formulation in both subdomains. Komatitsch \emph{et al.}~\cite{komatisch-tromp} introduced a Spectral Element method for modeling wave propagation in media with both fluid (acoustic) and solid (elastic) regions. The employed formulation is \emph{symmetric} (i.e., it is made in terms of displacement in elastic regions and velocity potential in acoustic regions), and matching between domains is implemented based on an interface integral in the framework of an explicit prediction-multicorrection staggered time scheme. Berm\'udez \emph{et al.}~\cite{bermudez2} considered a Finite Element approach to the problem based on a pressure formulation in the acoustic domain. Chaljub \emph{et al.}~\cite{capdeville} studied a Spectral Element approach for modeling elastic wave propagation in a solid-fluid sphere by taking into account the local effects of gravity, employing a symmetric formulation, as here. Flemisch \emph{et al.}~\cite{flemisch-wohlmuth} devised a numerical treatment based on two independent triangulations on the elastic and acoustic domains with Finite Elements. Due to the flexible construction of both grids, the finite element nodes on the elastic and acoustic boundary on the interface may, in general, not coincide, so as to allow as much flexibility as possible; as a result, non-conforming grids appear at the interface of the two subdomains. Käser and Dumbser \cite{kaser-dumbser} considered a numerical scheme suited for unstructured 2D and 3D meshes based on a Discontinuous Galerkin approach to simulate seismic wave propagation in heterogeneous media containing fluid-solid interfaces, using a formulation in terms of a first-order hyperbolic system in velocity-stress unknowns. The solution across element interfaces is handled by Riemann solvers or numerical fluxes. De Basabe and Sen \cite{debasabe-sen-2010} investigated the stability of the Spectral Element method and the Interior Penalty Discontinuous Galerkin method, considering the Lax-Wendroff method for time stepping and showing that it allows for a larger time step than the usual leap-frog finite difference method, with higher-order accuracy. Wilcox \emph{et al.}~\cite{ghattas} studied a high-order Discontinuous Galerkin scheme for the three-dimensional problem based on a velocity-strain formulation, allowing for the solution of the acoustic and elastic wave equations within the same unified framework, based on a first-order system of hyperbolic equations. Soares \cite{soares} considered a stabilized time-domain Boundary Element method to discretize each sub-domain. Bottero \emph{et al.}~\cite{bottero-komatitsch} used a time-domain Spectral Element method for simulations of wave propagation in the framework of ocean acoustics. Terrana \emph{et al.}~\cite{terrana} studied a high-order hybridizable Discontinuous Galerkin Spectral Element method, again based on a first-order hyperbolic velocity-strain formulation of the wave equations written in conservative form. Very recently, {Appel\"{o}} and Wang \cite{appelo-wang} devised an energy-based Discontinuous Galerkin approach, again using a symmetric formulation. Finally, a detailed $hp$-convergence analysis of a Discontinuous Galerkin method on polytopal meshes has been presented and validated in a two-dimensional setting in \cite{bonaldi-mox}, wherein also a well-posedness result has been obtained by a semigroup-based approach.

In this paper, the unknowns of the problem are the displacement field in the solid domain and the velocity potential in the fluid domain, i.e., we employ a symmetric formulation. The latter, say $\phi$, is defined in terms of the acoustic velocity field $\bf v_a$ in such a way that $\bf v_a = -\grad\phi$. Also, the acoustic pressure $p_a$ in the fluid region is given by $p_a =\rho_a \dot\phi$, with $\dot\phi$ the first time derivative of the velocity potential. 

In the context of earthquake ground motion simulations, the numerical scheme employed has to satisfy the following requirements: \emph{accuracy}, \emph{geometric flexibility}, and \emph{scalability}. To be \emph{accurate}, the numerical method must keep dissipative and dispersive errors low. \emph{Geometric flexibility} is required since the computational domain usually features complicated geometrical shapes as well as sharp discontinuities of mechanical properties. Finally, real-life seismic scenarios are typically characterized by domains whose dimension, ranging from hundreds to thousands square kilometers, is very large compared with the wavelengths of interest. This typically leads to a discrete problem featuring several millions of unknowns. As a consequence, parallel algorithms must be \emph{scalable} in order to efficiently exploit high performance computers.

To comply with these requirements, we employ a Discontinuous Galerkin Spectral Element (DGSE) approach based on a domain decomposition paradigm, which was introduced in \cite{rapetti}. More precisely, the discontinuities are imposed only at the interfaces between suitable non-conforming macroregions, so that the flexibility of the DG methods is preserved while keeping the accuracy and efficiency of Spectral Element (SE) methods and avoiding the proliferation of degrees of freedom that characterize DG methods. We refer to \cite{mox-elastodynamics} for a more detailed and comprehensive review of discretization methods for seismic wave propagation problems. 

The rest of the paper is organized as follows. In Section~\ref{sec:problem.statement} we give the formulation of the problem and recall the well-posedness result proven in \cite{bonaldi-mox} under suitable hypotheses on source terms and initial values. In Section~\ref{sec:disc.set} we introduce the DGSE method and present the formulation of the {semi-discrete} problem, also recalling a stability result for its formulation in a suitable energy norm, as well as $hp$-convergence results (with $h$ and $p$ denoting the meshsize and the polynomial approximation degree, respectively) for the error in the same norm; a discussion of the fully discrete formulation of the problem is presented as well. Finally, in Section~\ref{sec:numerical.results}, we present several numerical experiments carried out in a three-dimensional setting, with the two-fold aim of {verifying} the theoretical results and simulating test cases of physical interest. 

{Throughout the paper, we will use standard notation for Sobolev spaces \cite{adams-fournier}. The Sobolev spaces of vector-valued functions are denoted by $\textbf{H}^m(\Omega) \equiv [H^m(\Omega)]^d$ and their norms by $\|{\cdot}\|_{m,\Omega}$, where $\Omega\subset\mathbb{R}^d$ is an open bounded domain of $\bb R^d$, $d\in\{2,3\}$. We will use the symbol $(\cdot , \cdot)_\Omega$ and $\|{\cdot}\|_\Omega$ to denote the standard inner product and norm in the space $H^0 (\Omega) \equiv L^2 (\Omega)$, respectively.} We also use the abridged notation $x \lesssim  y$ in place of $x \le Cy$, for $C > 0$ independent of the discretization parameters (polynomial degree and meshsize), but possibly depending on the material properties of the media under consideration.

\section{Problem statement}\label{sec:problem.statement}
In this section, we recall the formulation of the elasto-acoustic problem in its \emph{symmetric} form, i.e., written in terms of the displacement field $\bf u$ and the velocity potential $\phi$, defined such that the velocity field in the acoustic domain $\bf v_a$ is given by $\bf v_a = -\grad\phi$ (see \cite{bonaldi-mox}).
Let $\Omega \equiv \Oe \cup \Oa \subset \bb R^3$ denote an open bounded domain
with Lipschitz boundary, given by the union of two open disjoint bounded subdomains $\Omega_e$ and $\Omega_a$ representing the elastic and acoustic regions in their reference configurations, respectively. We denote by $\Gamma_{\mr I} \equiv \de\Omega_e \cap \de\Omega_a$ the \emph{interface}
between the two domains. Thus, given a body force $\bf f_e$ and a scalar volume acoustic source $f_a$ as well as a final time $T>0$, the strong formulation of the problem reads
\begin{align}\label{strong_form}
\left\{
\begin{alignedat}{2}
{\rho_e \ddot{\bf u}  - \divb\bm\sigma(\bf u)} &= {\bf f_e} &\qquad&\hbox{in } \Omega_e \times(0,T],\\[4pt]
%
\bm\sigma(\bf u)\bf n_e & = - \rho_a\dot\phi\, \bf n_e &\qquad& \hbox{on }\Gamma_{\mr I}\times (0,T], \\[4pt]
c^{-2}\ddot\phi - \Del\phi &= f_a &\qquad&\hbox{in } \Omega_a \times (0,T],\\[4pt]
%
%
{\de\phi}/{\de\bf n_a} &= - \dot{\bf u} \dotp \bf n_a &\qquad& \hbox{on }\Gamma_{\mr I}\times (0,T],
 \\[4pt]
\end{alignedat}
\right.
\end{align}
{coupled with suitable boundary and initial conditions that are detailed below.} 

Here, $\rho_e \in L^\infty(\Oe)$, $\rho_e > 0$, is the {mass} density of the elastic region $\Oe$; $\bm\sigma(\bf u) = \bb C\bm\eps(\bf u) = \lambda(\tr\bm\eps(\bf u))\bf I + 2\mu\bm\eps(\bf u)$ is the {Cauchy stress} tensor; $\bb C$ is the uniformly elliptic and symmetric fourth-order elasticity tensor, representing a linearly elastic isotropic behavior, with $\mu$ and $\lambda$ the Lam\'e parameters; $\bm\eps(\bf u) = \sym(\grad\bf u) =  \frac{1}{2}\left(\grad \bf u + \grad \bf u^T\right)$ is the strain tensor. Also, we denote by $\rho_a\in L^\infty(\Oa)$, $\rho_a > 0$, the density of the acoustic region $\Oa$ and by $c > 0$ the speed of the acoustic wave. 

The trasmission conditions on $\Gamma_{\mr I}$ take account of the pressure, of magnitude $\rho_a |\dot\phi|$, exterted by the acoustic region onto the elastic one through the interface, and of the continuity of the normal component of the velocity field at the interface.

Concerning boundary conditions we assume the following decomposition: $\de \Omega = (\de \Omega_e \cup \de \Omega_a) \setminus \Gamma_{\mr I}$ where
$\de\Omega_e  = \Gamma_{e,D} \cup \Gamma_{e,N} \cup \Gamma_{e,N\!R} \cup \Gamma_{\mr I}$ and 
$\de\Omega_a  = \Gamma_{a,D} \cup \Gamma_{a,N} \cup \Gamma_{a,N\!R} \cup \Gamma_{\mr I}$. We denote by $\bf n_e$ and $\bf n_a$ the outer unit normal vectors to $\de\Omega_e$ and $\de\Omega_a$, respectively.
Homogeneous Dirichlet boundary conditions are assigned on $\Gamma_{e,D} \cup \Gamma_{a,D}$, i.e., $\bf u = \bf 0$ and $\phi = 0$.   Neumann boundary conditions on $\Gamma_{e,N} \cup \Gamma_{a,N}$ are prescribed in term of a surface traction  $\bf g_e$ and a surface acoustic flux $g_a$ as
\begin{align*}
\left\{
\begin{alignedat}{2}
\bm\sigma(\bf u)\bf n_e & = \bf g_e &\qquad& \hbox{on }\Gamma_{e,N} \times (0,T],\\
{\de\phi}/{\de\bf n_a} &= g_a &\qquad& \hbox{on }\Gamma_{a,N}\times (0,T]. \\
\end{alignedat}
\right.
\end{align*}
Non-reflecting boundary conditions are imposed on $\Gamma_{e,N\!R} \cup \Gamma_{a,N\!R}$; here, the surface loads are themselves expressed in terms of the time derivatives of the unknowns. In particular, we set 
\begin{equation}\label{abso}
\left\{
\begin{alignedat}{2}
\bm\sigma(\bf u)\bf n_e & = \bf g_e^\star &\qquad& \hbox{on }\Gamma_{e,\NR} \times (0,T],\\
{\de\phi}/{\de\bf n_a} &= g_a^\star &\qquad& \hbox{on }\Gamma_{a,\NR}\times (0,T], \\
\end{alignedat}
\right.
\end{equation}
with $\bf g_e^\star \equiv \rho_e c_P  (\dot{\bf u} \cdot\bf n_e)\bf n_e + \rho_e c_S \dot{\bf u}_\tau$ and $g_a^\star \equiv c^{-1}\dot\phi$
(see e.g.~\cite{speed, komatisch-vilotte, flemisch-wohlmuth}), with $\dot{\bf u}_\tau \equiv \dot{\bf u} - (\dot{\bf u} \cdot\bf n_e)\bf n_e$ is the tangential velocity field over $\Gamma_{e,N\!R}$, and $c_P$ and $c_S$ are the propagation velocities of P (pressure) and S (shear) waves, respectively, given by $c_P = \sqrt{{(\lambda+2\mu)}/{\rho_e}}$ and $ c_S = \sqrt{{\mu}/{\rho_e}}$.
These are commonly referred to in literature as {first order absorbing boundaries \cite{stacey}}. 

Finally, as initial conditions we set $\bf u(\cdot,0) = \bf u_0$ and $\dot{\bf u}(\cdot,0)  = \bf u_1$ in $ \Omega_e$ while 
$\phi(\cdot,0) = \phi_0$ and $\dot{\phi}(\cdot,0) = \phi_1 $  in $\Omega_a$, for some regular enough functions ${\bf u_0}, {\bf u_1}, \phi_0$, and $\phi_1$.


The well-posedness of the problem {\eqref{strong_form}} in suitable functional spaces was proven in \cite{bonaldi-mox} {under suitable regularity assumptions on the data}, in the case $\Gamma_{e,N}\cup\Gamma_{a,N}=\emptyset = \Gamma_{e,N\!R} \cup \Gamma_{a,N\!R}$.

\section{Numerical discretization}\label{sec:disc.set}
In this section we present the numerical approximation of the weak formulation of \eqref{strong_form} through a DGSE method coupled with an explicit \emph{Newmark predictor-{corrector} staggered} time marching scheme (see \cite{komatisch-tromp}).
We first introduce the semi-discrete counterpart of \eqref{strong_form}, observing that the solution of \eqref{strong_form} satisfies the following {weak form}: for any $t\in (0,T]$, and all $(\bf v, \psi) \in  {\bf H^1_{\Gamma_{e,D}}(\Oe)} \times {H^1_{\Gamma_{a,D}}(\Oa)}$,
\begin{equation}
\begin{multlined}\label{weak_form}
(\rho_e\ddot{\bf u}(t),\bf v)_{\Oe} +(c^{-2}\rho_a\ddot\phi(t),\psi)_{\Oa} +\cc A_e(\bf u(t),\bf v)  + \cc A_a(\phi(t),\psi) +\cc I_e(\dot\phi(t), \bf v) + \cc I_a(\dot{\bf u}(t), \psi) \\[2pt]  = ({\bf f_e}(t), \bf v)_{\Oe} + (\bf g_e(t), \bf v)_{\Gamma_{e,N}} + (\bf g_e^\star(t), \bf v)_{\Gamma_{e,N\!R}}  \\ + (f_a(t),\psi)_{\Oa} + (g_a(t),\psi)_{\Gamma_{a,N}} + (\bf g_a^\star(t), \bf v)_{\Gamma_{a,N\!R}},
\end{multlined}
\end{equation}
where
\begin{align}
\begin{aligned}
\cc A_e(\bf u, \bf v) & = ({\bb C}\bm\eps(\bf u),\bm\eps(\bf v))_{\Oe}, \\
\cc A_a(\phi,\psi) & = (\rho_a\grad\phi,\grad\psi)_{\Oa},
\end{aligned}
&\qquad
\begin{aligned} 
{\cc I_e(\psi,\bf v)}  & = (\rho_a\psi \bf n_e, \bf v)_{\Gamma_{\mr I}}, \\
{\cc I_a(\bf v, \psi)} & = (\rho_a {\bf v}\dotp \bf n_a,\psi)_{\Gamma_{\mr I}}.
\end{aligned}
\end{align}
We observe that the second evolution equation has been multiplied by $\rho_a$ to ensure (skew)symmetry of the two interface terms (since $\bf n_a = -\bf n_e$).
\subsection{Partitions and trace operators}
We now consider a decomposition $\cc T_{\Omega_e}$ of $\Oe$ into
 $L_e$ nonoverlapping polyhedral regions $\Omega_e^\ell$, $\ell \in \{1,\dots,L_e\}$, such that
$\Oe = \bigcup_{\ell=1}^{L_e} \Oe^\ell$, with $\Omega_e^\ell \cap \Omega_e^{\ell'} = \emptyset$ for any $\ell \neq \ell'$. This
first macropartition is introduced to distinguish elastic materials with different properties (density $\rho_e$ and material moduli $\lambda,\mu$).
On each $\Oe^\ell$, we build a conforming computational mesh $\cc T_{h_\ell}^e$ of meshsize $h_{\ell} > 0$
made of disjoint elements $\cc K_{e,\ell}$, and suppose that each $\cc K_{e,\ell} \subset \Oe^\ell$ is the image
through an invertible bilinear map $F_{e,\ell} : \wh{\cc K} \to \cc K_{e,\ell}$ of {the unit reference hexahedron
$\wh{\cc K} = (-1,1)^d$, $d\in\{2,3\}$}. Given two adjacent regions $\Oe^{\ell^\pm}$, we define an internal
face $F$ as the non-empty interior of $\de\wl{\K_{e,\ell^+}} \cap \de\wl{\K_{e,\ell^-}}$,
for  $\K_{e,\ell^\pm} \in {\cc T}_{h_{\ell^\pm}}^e$, $\K_{e,\ell^\pm} \subset \Omega_e^{\ell^\pm}$, and collect all the internal faces in
the set $\cc F_{h,e}^i$. Moreover, we define $\cc F_{h,e}^D$, $\cc F_{h,e}^N$, and $\cc F_{h,e}^\NR$ as the sets of all boundary faces
where displacement, tractions, or non-reflecting elastic boundary conditions are imposed, respectively. We collect all the boundary faces \emph{not} laying on $\Gamma_{\mr I}$ in the set $\cc F_{h,e}^b$.


On the other hand, concerning the acoustic domain $\Oa$, since we do not take into account multi-phase fluids, we introduce a
conforming grid $\cc T_h^a$ of $\Oa$ made by disjoint hexahedral elements $\cc K_a$. As in the elastic case, we suppose that each $\cc K_{e,\ell} \subset \Oe^\ell$ is the image through an invertible bilinear map $F_{a} : \wh{\cc K} \to \cc K_{a}$ of {the unit reference hexahedron
$\wh{\cc K} = (-1,1)^d$, $d\in\{2,3\}$}. 
Also, we define $\cc F_{h,a}^D$, $\cc F_{h,a}^N$, and $\cc F_{h,a}^\NR$ as the sets of all boundary faces
where velocity potential, fluxes, or non-reflecting acoustic boundary conditions are imposed, respectively.
We collect all the boundary faces \emph{not} laying on $\Gamma_{\mr I}$ in the set $\cc F_{h,a}^b$.

Finally, we collect all faces laying on $\Gamma_{\mr I}$ in the set $\cc F_{h,\Gamma_{\mr I}}$; {in this case, $F \in \cc F_{h,\Gamma_{\mr I}}$ is the non-empty interior of $\de\wl{\K_{e,\ell}} \cap \de\wl{\K_a}$, for given $\K_{e,\ell} \subset \Omega_e^\ell \in \cc T_{\Omega_e}$ and $\K_a \in \cc T_h^a$}. Implicit in these definitions is the assumption that each face laying on $\de\Oe \cup \de\Oa$ can belong to exactly one of the sets $\cc F_{h,e}^D$, $\cc F_{h,e}^N$, $\cc F_{h,a}^D$, $\cc F_{h,a}^N$, and $\cc F_{h,\Gamma_{\mr I}}$.

\begin{figure}
\includegraphics[scale=1.1]{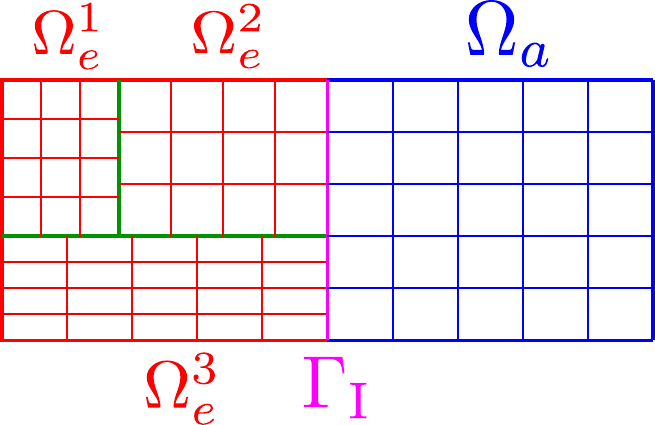}
\caption{Example of decompositions for the domains $\Oe$ and $\Oa$. Interfaces between elastic regions $\Oe^1$, $\Oe^2$, and $\Oe^3$, characterized by different material properties, are highlighted in green, and the elasto-acoustic interface $\Gamma_{\mr I}$ is highlighted in purple.}
\label{griglia}
\end{figure}

\begin{remark}[Non-matching grids at the elasto-acoustic interface]
\label{interface}
Notice that the above-detailed framework allows to handle the situation of non-matching grids at the interface $\Gamma_\mr I$ between the elastic and the acoustic domains (cf.~Figure~\ref{griglia}). Meshes can therefore be generated independently on each of the domains.
\end{remark}

We now introduce the following average and jump operators \cite{unified-brezzi,arnold-RM} {for mesh faces in the elastic domain}. For sufficiently smooth scalar, vector, and tensor fields $\psi$, $\bf v$, and $\bm\tau$, we define averages and jumps on an internal face{ $F \in \cc F_{h,e}^i$, 
$ F \subset \de \cc K_{e,\ell^+} \cap \de \cc K_{e,\ell^-} $}, with $\cc K_{e,\ell^\pm} \in \cc T_{h_{\ell^\pm}}^e$, as follows:
\begin{equation*}
\begin{alignedat}{2}
\jm{\psi} &= \psi^+\bf n^+ + \psi^- \bf n^-,&\qquad \av{\psi} &= \frac{\psi^+ + \psi^-}{2},\\
\jm{\bf v} & = \bf v^+ \otimes \bf n^+ + \bf v^- \otimes \bf n^-,&\qquad \av{\bf v} & = \frac{\bf v^+ + \bf v^-}{2},\\
\jm{\bm\tau} & = \bm\tau^+\bf n^+ + \bm\tau^-\bf n^-,&\qquad \av{\bm\tau} & = \frac{\bm\tau^+ + \bm\tau^-}{2},
\end{alignedat}
\end{equation*}
where $\bf a\otimes \bf b$ denotes the tensor product of $\bf a,\bf b \in \bb R^3$; $\psi^\pm$, $\bf v^\pm$ and $\bm\tau^\pm$ are the traces of $\psi$, $\bf v$ and $\bm\tau$ on $F$ taken from the interior of {$\cc K_{e,\ell^\pm}$}, and $\bf n^\pm$ is the outer unit normal vector to {$\de\cc K_{e,\ell^\pm}$}. When considering a {boundary face} {$F \in \cc F_{h,e}^b$}, we set $\jm{\psi} = \psi\bf n$, $\jm{\bf v} = \bf v \otimes \bf n$, $\jm{\bm\tau} = \bm\tau\bf n$, and $\av{\psi} = \psi$, $\av{\bf v} = \bf v$, $\av{\bm\tau} = \bm\tau$. We also use the shorthand notation 
$$\la  \Phi, \Psi\ra_{\cc F} = \sum_{F\in\cc F} (\Phi, \Psi)_F,
\qquad 
\|\Phi\|_{\cc F} = \la \Phi, \Phi\ra_{\cc F}^{\nf12}, $$ 
for scalar, vector or tensor fields $\Phi$ and $\Psi$ and for a given generic collection $\cc F$ of mesh faces.

\subsection{Discontinuous Galerkin Spectral Element approximation}
First, concerning the elastic domain $\Oe$, we associate with each subdomain $\Oe^\ell$ 
a nonnegative integer ${N_{e,\ell} \ge 1}$, and introduce the finite-dimensional space {
\begin{equation}
\label{V_e}
\bf V(\Oe^\ell) = \{ \bf v \in \bf C^0(\wl{\Omega}_e^\ell) : \bf v_{|\cc K_{e,\ell}} \circ F_{e,\ell} \in [\bb Q^{{N_{e,\ell}}}(\wh{\cc K})]^d  \  \forall \cc K_{e,\ell} \in {\cc T}_{h_\ell}^e \},
\end{equation}
where $\bb Q^{{N_{e,\ell}}}(\wh{\cc K})$ is  the space  of polynomials of degree $N_{e,\ell}$ in each coordinate direction on 
the unit reference hexahedron $\wh{\cc K}$}. We then introduce the space {$\bf V(\Oe)
= \bigtimes_{\ell=1}^{L_e} \bf V(\Oe^\ell)$}. Concerning the acoustic domain $\Oa$ we choose a spectral degree $N_a \ge 1$ and define the following space:
\begin{equation}\label{V_a}
{
V(\Oa) = \{ \psi \in C^0(\wl{\Omega}_a) : \psi_{|\cc K_a} \circ F_{a} \in \bb Q^{{N_a}}(\wh{\cc K})  \  \forall \cc K_a \in {\cc T}_{h}^a \}.
}
\end{equation}

The semi-discrete DGSE approximation of \eqref{weak_form} reads then: 
 $(\bf u_h, \phi_h)\! \in \! C^2([0,T];\! { \bf V(\Oe)})  \times C^2([0,T]; \! { V(\Oa))}$ such that, for all $(\bf v_h,\psi_h) \in  {\bf V(\Oe)} \times {V(\Oa)}$,
\begin{equation}\label{coupled_dG}
\begin{multlined}
(\rho_e\ddot{\bf u}_h(t),\bf v_h)_{\Oe} + (c^{-2}\rho_a\ddot\phi_h(t),\psi_h)_{\Oa}   
 + \cc A_h^e(\bf u_h(t),\bf v_h)  + \cc A_h^a(\phi_h(t),\psi_h) \\[5pt]
+  \cc I_h^e(\dot\phi_h(t), \bf v_h) + \cc I_h^a(\dot{\bf u}_h(t), \psi_h)
   = \cc L^e_h(\bf v_h) + \cc L^a_h(\psi_h), 
\end{multlined}
\end{equation}
with initial conditions $\left(\bf u_h (0),\dot{\bf u}_h(0)\right) = \left (\bf u_{0,h}, \bf u_{1,h}\right) \in { {\bf V}(\Oe)} \times {{\bf V}(\Oe)}$, and $\left( \phi_h(0),\dot{\phi}_h(0)\right) = \left(\phi_{0,h},\phi_{1,h} \right) \in  {{V}(\Oa)}\times {{V}(\Oa)}$,  where $\bf u_{0,h},\bf u_{1,h},\phi_{0,h},\text{ and }\phi_{1,h}$ are suitable approximations of the initial data. In \eqref{coupled_dG} 
\begin{equation}\label{discr_bilin_forms}
\begin{alignedat}{2}
 \cc A_h^e(\bf u,\bf v) & =  \sum_{\Omega_e^\ell \in \cc T_{\Omega_e}} (\bm\sigma_h(\bf u ),\bm\eps_h(\bf v))_{\Oe^\ell} - \la \av{\bm\sigma_h(\bf u)}, \jm{\bf v} \ra_{{\cc F}_{h,e}^i} \\  & \ \ \ - \la \jm{\bf u }, \av{\bm\sigma_h(\bf v)}\ra_{{\cc F}_{h,e}^i} + \la \eta\jm{\bf u },\jm{\bf v}\ra_{{\cc F}_{h,e}^i} &\ & \forall \bf u, \bf v \in {{\bf V}(\Oe)},\\[5pt]
 \cc A_h^a(\phi,\psi) & =  \sum_{\cc K_a \in \cc T_h^a} (\rho_a \grad\phi,\grad\psi)_{\cc K_a}  &\ &\forall \phi,\psi\in {V(\Oa)},\\[5pt]
 \cc I_h^e(\psi,\bf v) & =  \la\rho_a \psi \bf n_e,\bf v\ra_{\cc F_{h,\Gamma_{\mr I}}}&\ & \forall (\psi,\bf v) \in {V(\Oa)} \times {{\bf V}(\Oe)}, \\[5pt]
 \cc I_h^a(\bf v,\psi) & =  \la {\rho_a \bf v}\dotp{\bf n_a},\psi \ra_{\cc F_{h,\Gamma_{\mr I}}} \equiv -\cc I_h^e(\psi,\bf v)
 &\quad& \forall (\bf v,\psi) \in {{\bf V}(\Oe)} \times {V(\Oa)}, \\[2ex]
\cc L^e_h(\bf v) & =  \sum_{\Oe^\ell \in \cc T_{\Omega_e}} (\bf f_e(t) , \bf v)_{\Oe^\ell} + \la\bf g_e(t),\bf v \ra_{\cc F_{h,e}^N}
+ \la\bf g_e^\star(t),\bf v \ra_{\cc F_{h,e}^\NR}
&\quad&\forall \bf v \in {{\bf V}(\Oe)},\\
\cc L^a_h(\psi) & =  \sum_{\cc K_a \in \cc T_h^a} (f_a(t),\psi)_{\cc K_a} + \la g_a(t),\psi\ra_{\cc F_{h,a}^N} + \la g_a^\star(t),\psi\ra_{\cc F_{h,a}^\NR}
&\quad&\forall \psi \in  {V(\Oa)}.
\end{alignedat}
\end{equation}

We point out that the fourth identity in \eqref{discr_bilin_forms} holds since $\bf n_a = - \bf n_e$ on $\Gamma_{\mr I}$.
Here we have set, for any $\bf v \in {{\bf V}(\Oe)}$,
$$
\bm\eps_h(\bf v) = \frac{1}{2}\left(\grad_h \bf v + \grad_h \bf v^T\right),\quad
\bm\sigma_h(\bf v) = \bb C \bm\eps_h(\bf v),
$$
with $\grad_h$ the usual broken gradient operator.
The discontinuity penalization function $\eta\colon \cc F_{h,e}^i \to \bb R$ is defined as follows:
{
\begin{equation}\label{stabiliz.a}
\eta_{| F} =
 \alpha \{(\lambda+2\mu)^\pm\}_{\mr H} \left\{  \left(\frac{N_{e,\ell}^2}{h_\ell}\right)^\pm\right\}_{\mr H}, \quad 
F \in {\cc F_{h,e}^i},  \; F \subset \de{\cc K_{e,\ell^+}} \cap \de{\cc K_{e,\ell^-}},  
\end{equation}
with $\cc K_{e,\ell^\pm} \in \cc T_{h_{\ell^\pm}}^e$. Here,} $\alpha > 0$ is a positive constant to be properly chosen, and
$\{v^\pm \}_{\mr H} = 2 v^+ v^-/(v^+ + v^-)$ is the harmonic mean of traces $v^+$ and $v^-$ of a given scalar field $v$.

Upon introducing the following norms
\begin{equation}\label{norme}
\begin{alignedat}{2}
\|\bf v\|_{\DG,e}^2 & = \|\bb C^{\nf{1}{2}}\bm\eps_h(\bf v)\|_{\Oe}^2 + \|{\eta^{\nf{1}{2}}\jm{\bf v}}\|_{{\cc F_{h,e}^i}}^2
&\quad& \forall \bf v \in { {\bf V}(\Oe)},\\
\| \bf v(t)\|_{\E_e}^2 & = 
\| \rho_e^{\nf{1}{2}}\dot{\bf v}(t)\|_{\Oe}^2 + \|\rho_e^{\nf{1}{2}}\z \bf v(t)\|_{\Oe}^2 + \|\bf v(t)\|_{\DG,e}^2
&\quad& \forall \bf v \in C^1([0,T]; {\bf V(\Oe)}),\\[5pt] 
\| \psi(t)\|_{\E_a}^2 &=
\|c^{-1}\rho_a^{\nf{1}{2}}
\dot\psi(t)\|_{\Oa}^2 + \|\rho_a^{\nf12}\grad\psi(t)\|_{\Oa}^2 &\quad& \forall \psi \in C^1([0,T];{V(\Oa)}),
\end{alignedat}
\end{equation}
it is possible to prove that bilinear forms $\cc A_h^e$ and $\cc A_h^a$ are continuous and coercive. Consequently, a stability result and an error estimate in the above-defined energy norm for the semi-discrete solution can be inferred. We recall those results below; for the sake of readibility, we give a simplified statement of the error estimate (see~\cite{bonaldi-mox} for a more general framework, and \cite{ayuso} for the purely elastic case).
\begin{theorem}[Stability of the semi-discrete formulation]\label{stability}
Let $(\bf u_h,\phi_h)$ be the solution of \eqref{coupled_dG}. For a sufficiently large penalty parameter $\alpha$ in \eqref{stabiliz.a}, the following bound holds:
\begin{equation}\label{eq:stability}
\begin{aligned}
\|\bf u_h(t)\|_{\E_e} + \|\phi_h(t)\|_{\E_a} \lesssim  \|\bf u_h(0)\|_{\E_e} + \|\phi_h(0)\|_{\E_a} 
 + \int_0^t \left( \|{\bf f_e}(\tau)\|_{\Oe} + \|f_a(\tau)\|_{\Oa} \right)\,\d\tau,\quad t \in (0,T].
\end{aligned}
\end{equation}
\end{theorem}
\begin{theorem}[\emph{A priori} error estimate in the energy norm]\label{th-conve}
Assume that the exact solution of problem \eqref{strong_form} is such that
$\bf u \in C^2([0,T]; \bf H^m(\Oe))$ and $\phi \in C^2([0,T]; H^n(\Oa))$, for given integers $m,n \ge 2$.
Then, the following error estimate holds:
\begin{multline*}
\sup_{t\in[0,T]} (\|\bf u_h(t) - \bf u(t)\|_{\E_e}^2+\|\phi_h(t) - \phi(t)\|_{\E_a}^2)\\ \lesssim \sup_{t\in[0,T]}\left(\sum_{\ell =1}^{L_e} \frac{h_\ell^{2\min(m,N_{e,\ell}+1)-2}}{N_{e,\ell}^{2m -3}}\left(
\|\dot{\bf u}\|_{m,\Oe^\ell}^2 + \| \bf u\|_{m,\Oe^\ell}^2\right)
\!+\! \sum_{\K\in\cc T_h^a}\!\! \frac{h_\K^{2\min(n,N_a+1)-2}}{N_{a,}^{2n-3}}
\left( \| \dot\phi\|_{n,\cc K}^2 +  \| \phi\|_{n,\cc K}^2\right)
\right) 
\\
+ \int_0^T\!\left( \sum_{\ell=1}^{L_e} \frac{h_\ell^{2\min(m,N_{e,\ell}+1)-2}}{N_{e,\ell}^{2m-3}}\left(
\|\ddot{\bf u}\|_{m,\Oe^\ell}^2 +
 \|\dot{\bf u}\|_{m,\Oe^\ell}^2 + \| \bf u\|_{m,\Oe^\ell}^2\right)  \right. \\
\left. + \sum_{\K\in\cc T_h^a}\!\! \frac{h_\K^{2\min(n,N_{a,}+1)-2}}{N_{a}^{2n-3}}\left(
\|\ddot{\phi}\|_{n,\cc K}^2 +
\| \dot{\phi}\|_{n,\cc K}^2 + 
\| \phi\|_{n,\cc K}^2\right)\right) \d \tau. 
\end{multline*}
\end{theorem}
\begin{remark}[Error in the energy norm]
If both meshsizes are quasi-uniform, i.e. $h_\ell \simeq h_e\ \forall \ell \in \{1,\dots,L_e\}$ and $h_\K \simeq h_a\ \forall \cc K \in \cc T_h^a$, if the polynomial degree is uniform over elastic regions $\Omega_e^\ell$, i.e. $N_{e,\ell} = N_e\ \forall \ell \in \{1,\dots,L_e\}$, and if $m \ge N_e +1$ and $n \ge N_a+1$, the following error estimate holds:
\begin{equation}\label{errore'}
\begin{aligned}
\sup_{t\in[0,T]} (\|\bf u_h(t) - \bf u(t)\|_{\E_e} + \|\phi_h(t)-\phi(t)\|_{\E_a}) \lesssim \ {C_{\bf u}(T) \frac{h_e^{N_e}}{N_e^{m-\nf32}}}
+ {C_{\phi}(T) \frac{h_a^{N_a}}{N_a^{n-\nf32}}},
\end{aligned}
\end{equation}
where $C_{\bf u}(T)$ and $C_{\phi}(T)$ are positive numbers depending on the final time $T$ and the exact solution, along with its time derivatives.
\end{remark}
\subsection{Fully discrete formulation}\label{sec:fully.disc}
Upon fixing {polynomial} bases for discrete spaces ${\bf V(\Oe)}$ and ${V(\Oa)}$, { see e.g. \cite{bonaldi-mox}}, the semi-discrete algebraic formulation of problem \eqref{coupled_dG} reads
{
\begin{equation}\label{time_discr}
\begin{dcases}
\begin{aligned}
\sF M_e \ddot{\mathsf u}(t)  + \sF S_e\dot{\mathsf u}(t) +  \sF K_e 
\mathsf u(t) + \mathsf C_e \dot{\mathsf \upphi} (t)
& = \sF f_e(t), \quad t \in (0,T], \\
\mathsf M_a\ddot{\mathsf \upphi}(t) + \sF S_a\dot{\sF \upphi}(t) +
 \sF K_a\mathsf \upphi(t) + \mathsf C_a \dot{\mathsf u}(t) & = \sF f_a(t), \quad t \in (0,T],
\end{aligned}
\end{dcases}
\end{equation}
with initial conditions $\mathsf u(0)  = \mathsf u^0,\,
\dot{\mathsf u}(0) = \mathsf v^0, \, \mathsf{\upphi}(0)  = \mathsf  \upphi^0$, and $\dot{\mathsf\upphi}(0)  = \mathsf \uppsi^0$,} and where {the} vectors $\sF u(t)$ and $\sF \upphi(t)$ represent the expansion coefficients of $\bm u_h(t)$ and $\phi_h(t)$ in the chosen bases, respectively.
Analogously, $\sF M_e$, $\sF K_e$, and $\sF C_e$ are the matrix representations of the bilinear forms
{$(\rho_e\bf u,\bf v)_{\Oe},  \
\cc A_h(\bf u, \bf v)$ and  $\cc I_h^e(\psi,\bf v)$,
respectively (see~\eqref{discr_bilin_forms})}. When elastic absorbing boundary conditions are included in the model, matrix $\sF S_e$ 
takes account of the boundary term $(\bf g_e^\star,\bf v)_{\Gamma_{e,N\!R}}$; otherwise, it is identically equal to zero.
On the other hand, $\sF M_a$, $\sF K_a$, and $\sF C_a \equiv - \sF C_e^{\sF T}$ represent the bilinear forms {
$(c^{-2}\rho_a\phi,\psi)_{\Oa},\ \cc A_h^a(\phi,\psi)$, and $\cc I_h^a(\bf v,\psi)$,} respectively. When acoustic absorbing boundary conditions are considered, $\sF S_a$ represents the boundary term $(g_a^\star,\psi)_{\Gamma_{a,N\!R}}$. Finally, $\sF f_e(t)$ and $\sF f_a(t)$ are the vector representations of linear functionals $\cc L_h^e$ and $\cc L_h^a$, respectively.

For the time integration of system \eqref{time_discr}, as in \cite{komatisch-tromp}, we employ an explicit {Newmark predictor-{corrector} staggered} method \cite{hughes}; in this case, the scheme is conditionally stable and second-order accurate. 
We thus subdivide the time interval $[0,T]$ into $N_T$ subintervals of amplitude $\Delta t = T/N_T$ and denote by $\sF u^n \approx \sF u(t_n)$, $\sF v^n \approx \dot{\sF u}(t_n)$, $\sF a_e^n \approx \ddot{\sF u}(t_n)$, $\sF \upphi^n \approx \sF \upphi(t_n)$, $\sF{\uppsi}^n \approx \dot{\sF \Phi}(t_n)$, and $\sF a_a^n \approx \ddot{\sF \upphi}(t_n)$ the approximations of $\sF u$, $\dot{\sF u}$, $\ddot{\sF u}$, $\sF \upphi$, $\dot{\sF \upphi}$, and $\ddot{\sF \upphi}$ at time $t_n = n\Delta t$, $n \in \{0, \dots, N_T\}$, respectively. Then, along the lines of \cite{komatisch-tromp}, we exploit the fact that mass matrices are diagonal, and implement an iterative scheme based on a staggered {prediction/correction} technique. At each time step, we first compute predictors of the solution in both domains:
\begin{equation}\label{predictors}
\begin{alignedat}{3}
\wt{\sF u}^{n+1} & = \sF u^n + \Delta t \sF v^n + \frac{\Delta t^2}{2} \sF a_e^n, &\qquad& \wt{\sF v}^{n+1}  & = \sF v^n + \frac{\Delta t}{2} \sF a_e^n,\\
 \wt{\sF \upphi}^{n+1} & = \sF \upphi^n + \Delta t \sF \uppsi^n + \frac{\Delta t^2}{2} \sF a_a^n, &\qquad& \wt{\sF \uppsi}^{n+1} & = \sF \uppsi^n + \frac{\Delta t}{2} \sF a_a^n.
 \end{alignedat}
 \end{equation}
 Then, we update the solution in the elastic domain by solving the first equation of \eqref{time_discr} for $\sF a_e^{n+1}$, where the coupling term is evaluated as $-\sF C_e\wt{\uppsi}^{n+1}$, hence  {using the predictor computed in the acoustic domain}. Next, we compute the solution in the acoustic domain by solving the second equation of \eqref{time_discr} for $\sF a_a^{n+1}$, now using the updated solution in the elastic domain to evaluate the coupling term, which is thus given by $-\sF C_a \sF v^{n+1}$, where $\sF v^{n+1} = \wt{\sF v}^{n+1} + \frac{\Delta t}{2} \sF a_e^{n+1}$. We then iterate this algorithm by returning to the first step, this time using the updated solution.
 {The algorithm is summarized in the following scheme.}
 
 \begin{tcolorbox}
	\textbf{Newmark predictor-corrector staggered scheme}
	\\
	Given initial conditions $\sF u^0, \sF v^0$ and $\sF\upphi^0, \sF\uppsi^0$:
	\begin{algorithmic}
	\STATE compute $\sF a_e^0$ and $\sF a_a^0$:
	$\left\{
	\begin{aligned}
	\sF M_e \sF a_e^0 & = \sF f_e^0 - \sF S_e \sF v^0 - \sF K_e \sF u^0 - \sF C_e \sF \uppsi^0,\\
	\sF M_a \sF a_a^0 & = \sF f_a^0 - \sF S_a\sF \uppsi^0 - \sF K_a \sF \upphi^0 - \sF C_a \sF v^0;
	\end{aligned}
	\right.
	$ \\[1ex]
	\FOR{$n=0$ \TO $N_T-1$}
	\STATE compute predictors $\wt{\sF u}^{n+1},\wt{\sF v}^{n+1},\wt{\sF \upphi}^{n+1},\wt{\sF \uppsi}^{n+1}$ as in \eqref{predictors};
	\STATE compute $\sF a_e^{n+1}$:
	$\sF M_e \sF a_e^{n+1} = \sF f_e^{n+1} - \sF S_e \wt{\sF v}^{n+1} - \sF K_e \wt{\sF u}^{n+1} - \sF C_e \wt{\sF\uppsi}^{n+1};$
	\STATE update the solution in $\Omega_e$:
	$\sF u^{n+1}  = \wt{\sF u}^{n+1}, \ \sF v^{n+1}= \wt{\sF v}^{n+1} + \frac{\Delta t}{2}\sF a_e^{n+1}$;
	\STATE compute $\sF a_a^{n+1}$:
	$\sF M_a \sF a_a^{n+1} = \sF f_a^{n+1} - \sF S_a \wt{\sF \uppsi}^{n+1} - \sF K_a \wt{\sF \upphi}^{n+1} - \sF C_a {\sF v}^{n+1};$
	\STATE update the solution in $\Omega_a$: 
	$\sF \upphi^{n+1}  = \wt{\sF \upphi}^{n+1}, \ \sF \uppsi^{n+1}= \wt{\sF \uppsi}^{n+1} + \frac{\Delta t}{2}\sF a_a^{n+1}$;
	\ENDFOR
	\STATE
\end{algorithmic}
\end{tcolorbox}
 
 
\section{Numerical results}\label{sec:numerical.results}
\subsection{Verification {test}}
\label{sec_test}
\begin{figure}
\centering
\subfloat[$h=h_e=h_a=0.1$]{\includegraphics[keepaspectratio=true,scale=.275]{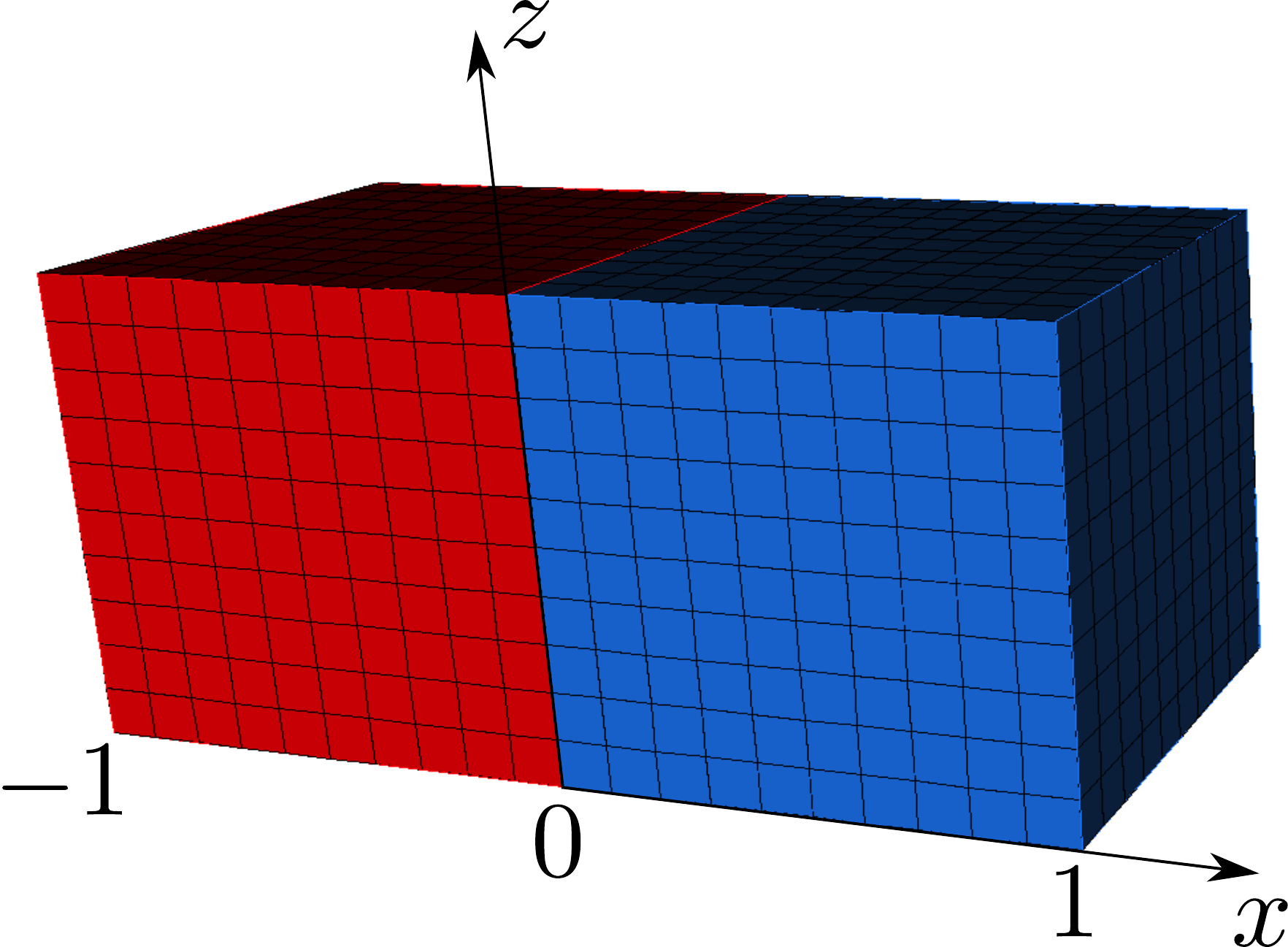}}
\subfloat[$h_e=0.1$, $h_a=0.2$]{\includegraphics[keepaspectratio=true,scale=.225]{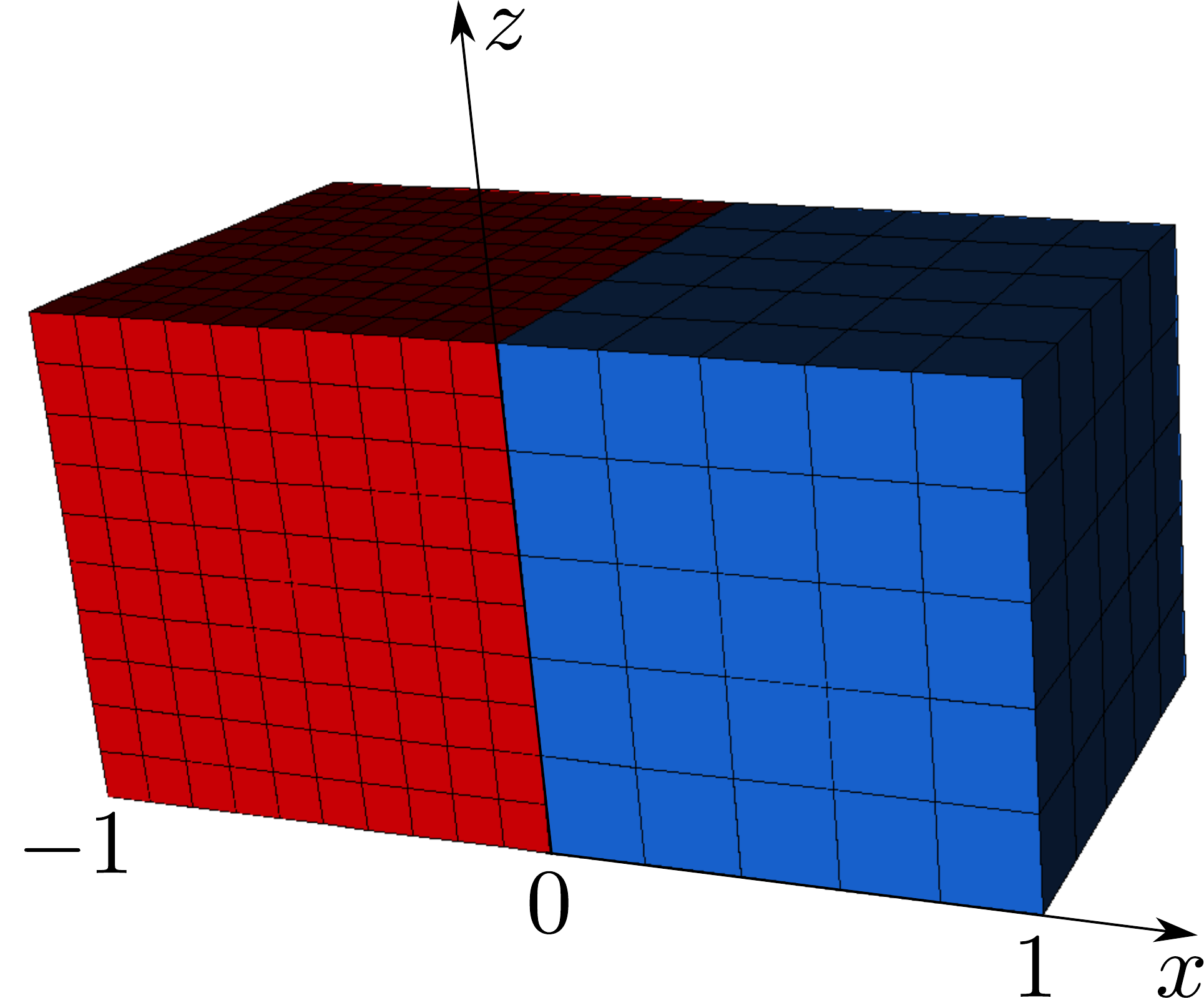}}
\subfloat[$h_e=0.1$, $h_a=0.15$]{\includegraphics[keepaspectratio=true,scale=.25]{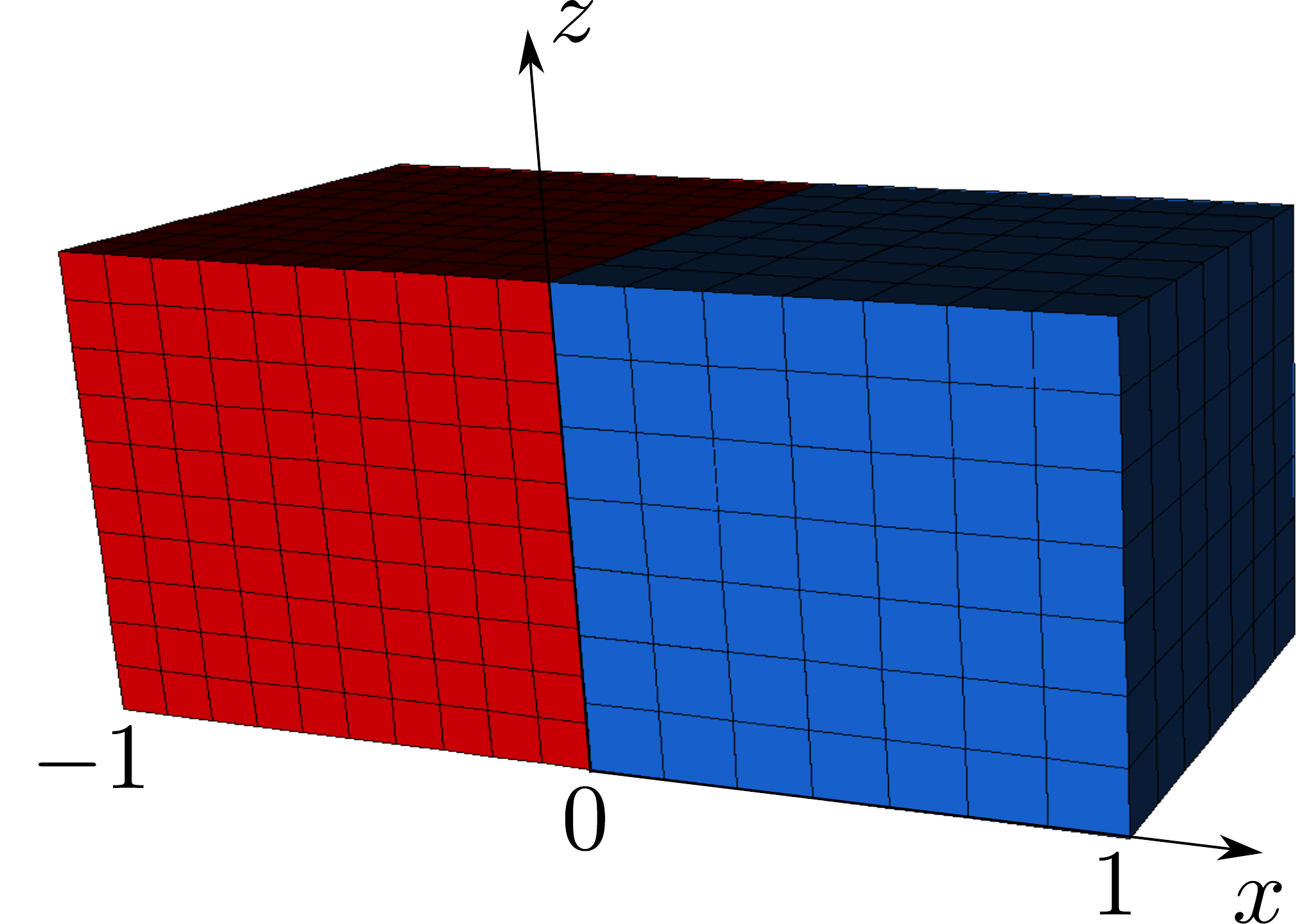}}
\caption{Test case~\ref{sec_test}. Computational domain with matching (a) and non-matching (b)--(c) hexahedral meshes.}
\label{mesh}
\end{figure}

\begin{figure}
\centering
\subfloat[$N=N_e=N_a=2$]{\includegraphics[scale=.475]{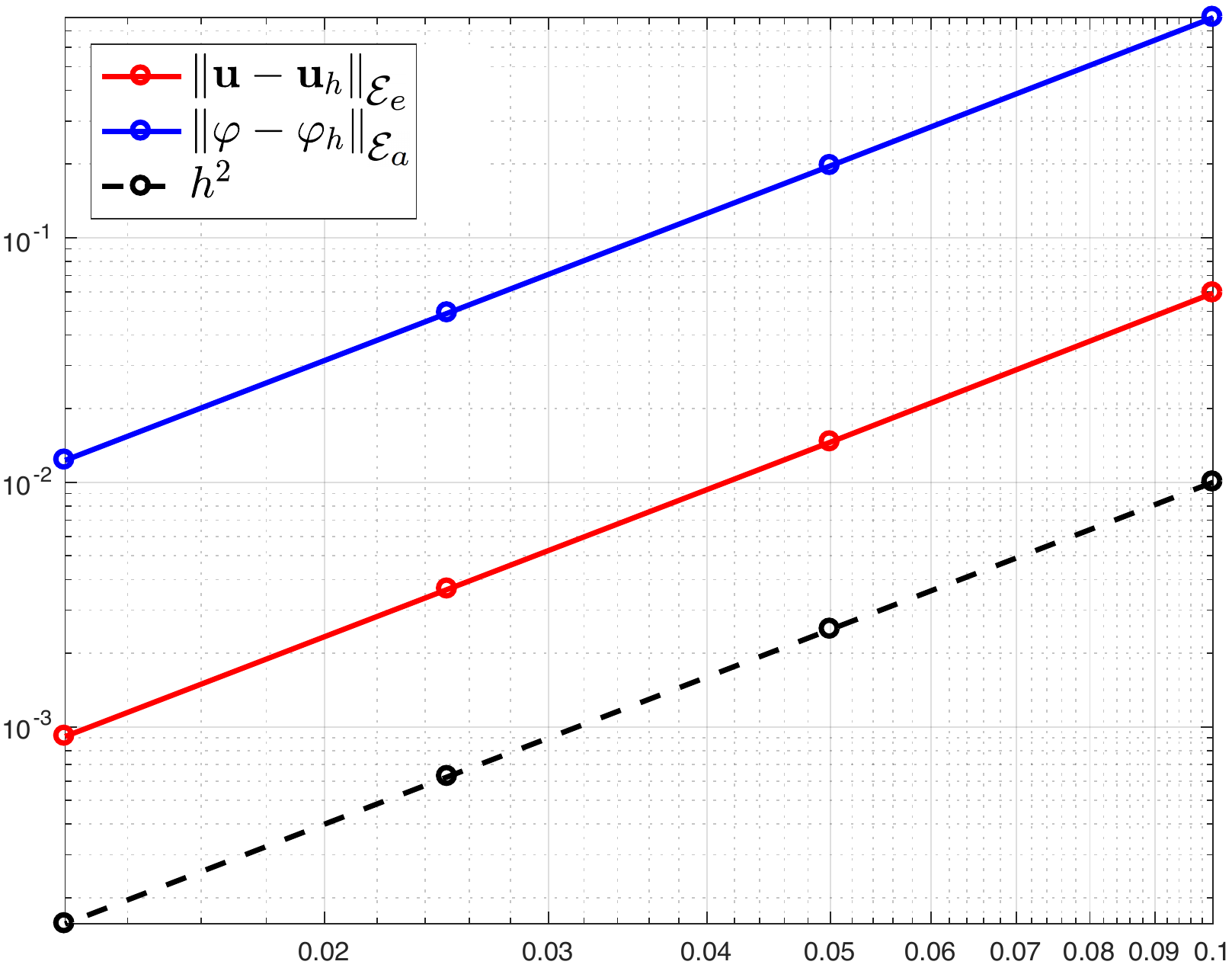}}
\subfloat[$N=N_e=N_a=3$]{\includegraphics[scale=.475]{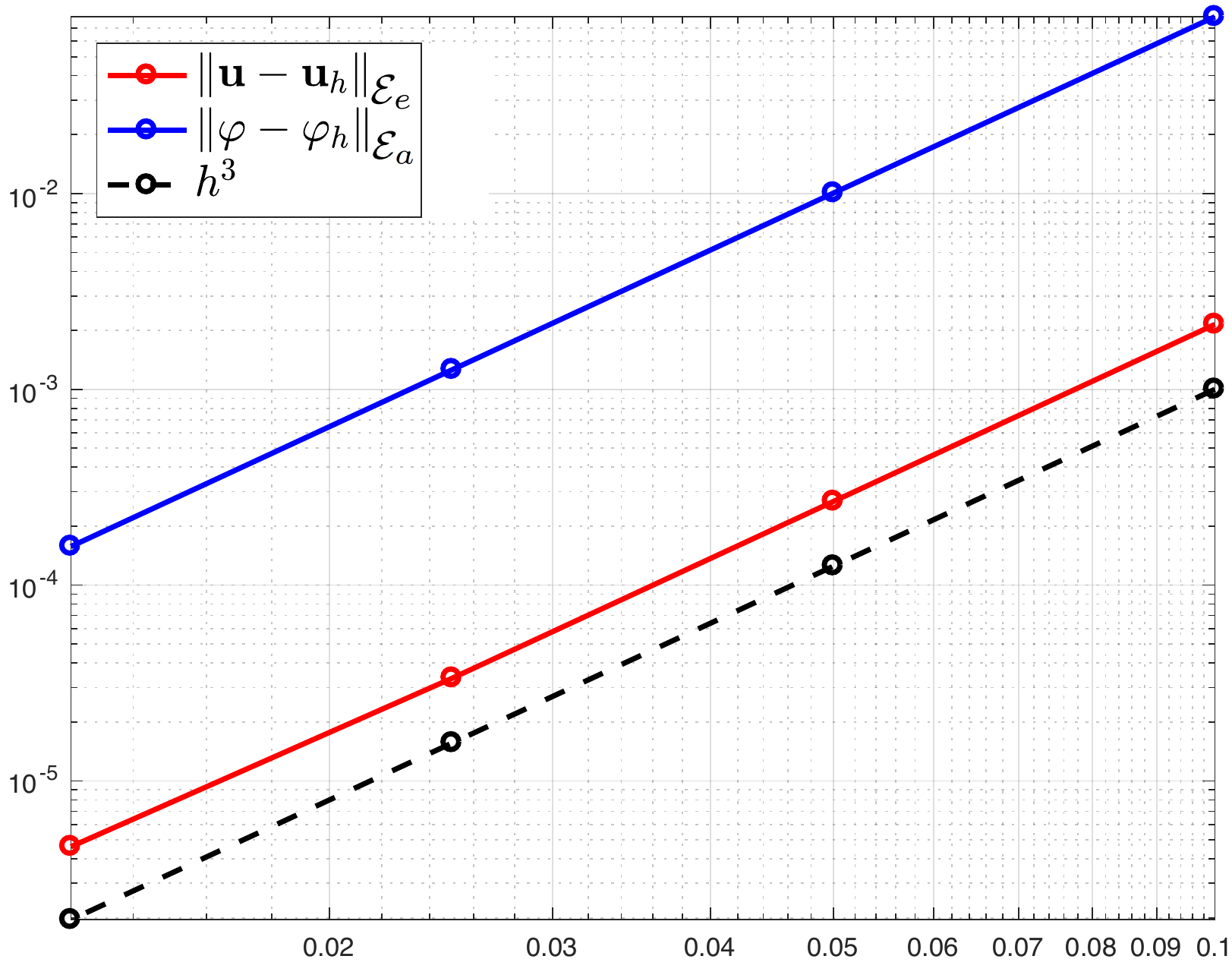}}\\
\subfloat[$N=N_e=N_a$ ranging from 2 to 7]{\includegraphics[scale=.475]{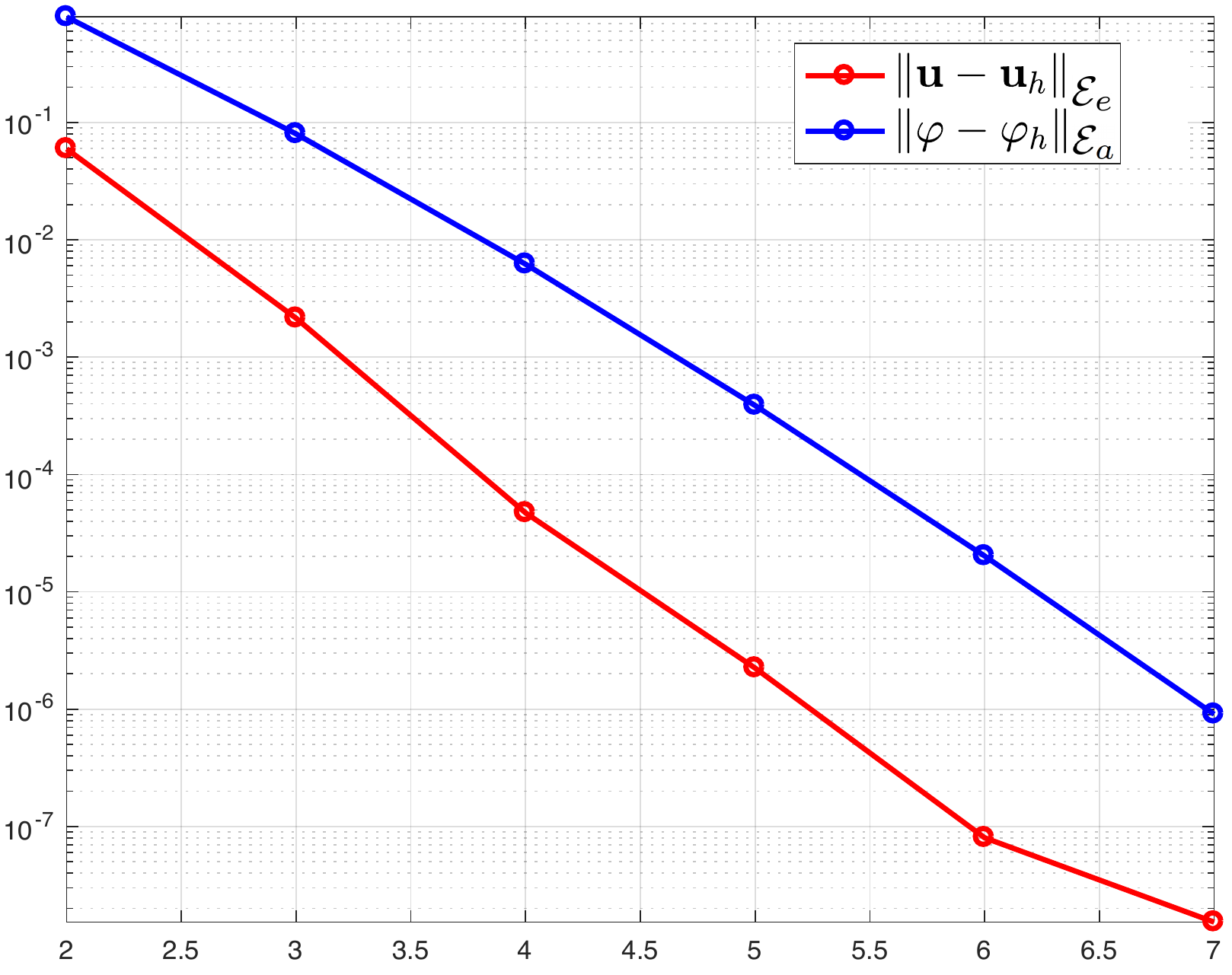}}
\caption{Test case~\ref{sec_test}. Error in the energy norm vs.~$h$ (a)--(b) and $N$ (c) at $t=0.1\,s$.}
\label{errore_h}
\end{figure}

\begin{figure}
\centering
\subfloat[$N=N_e=N_a=2$]{\includegraphics[scale=.475]{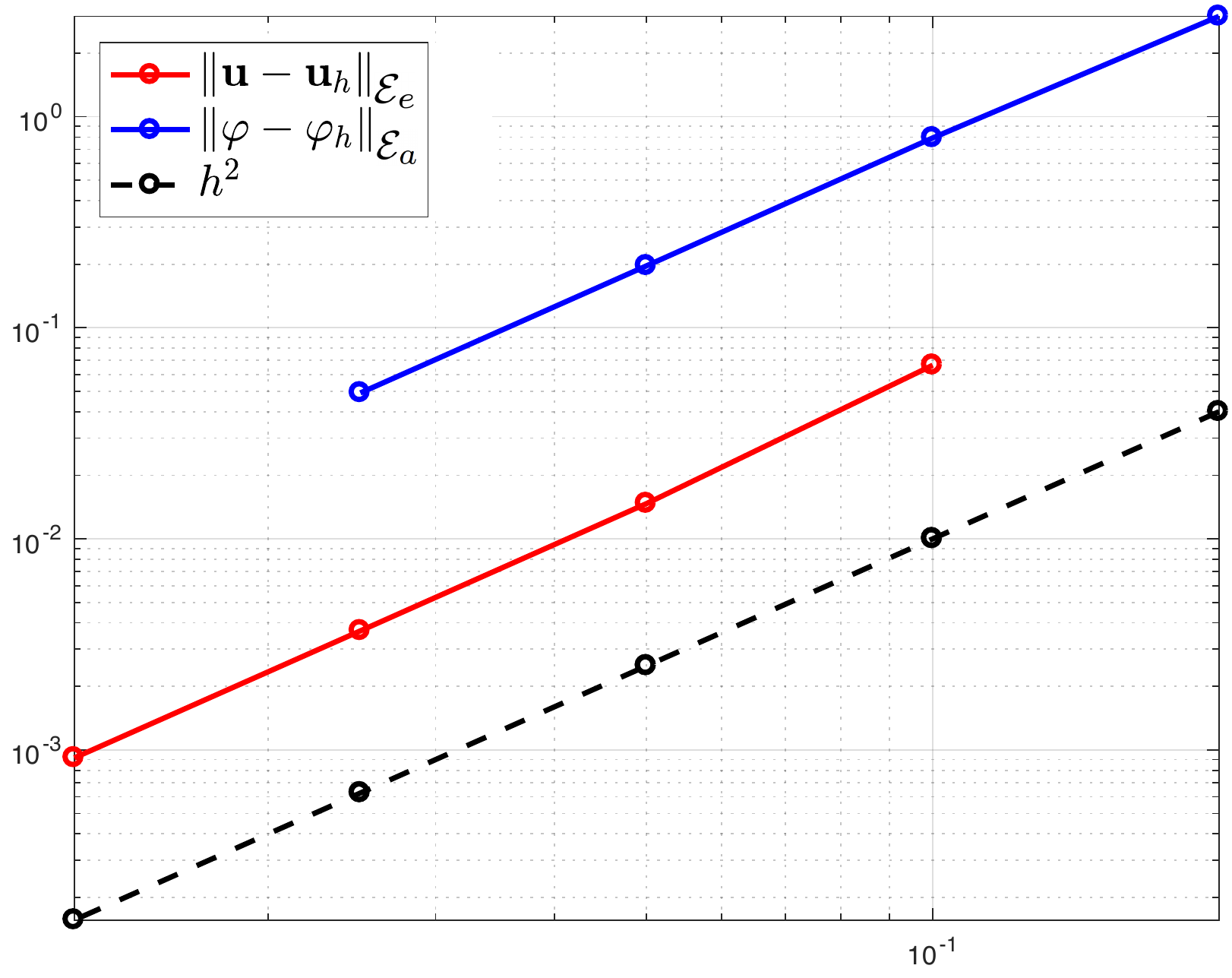}}
\subfloat[$N=N_e=N_a=3$]{\includegraphics[scale=.475]{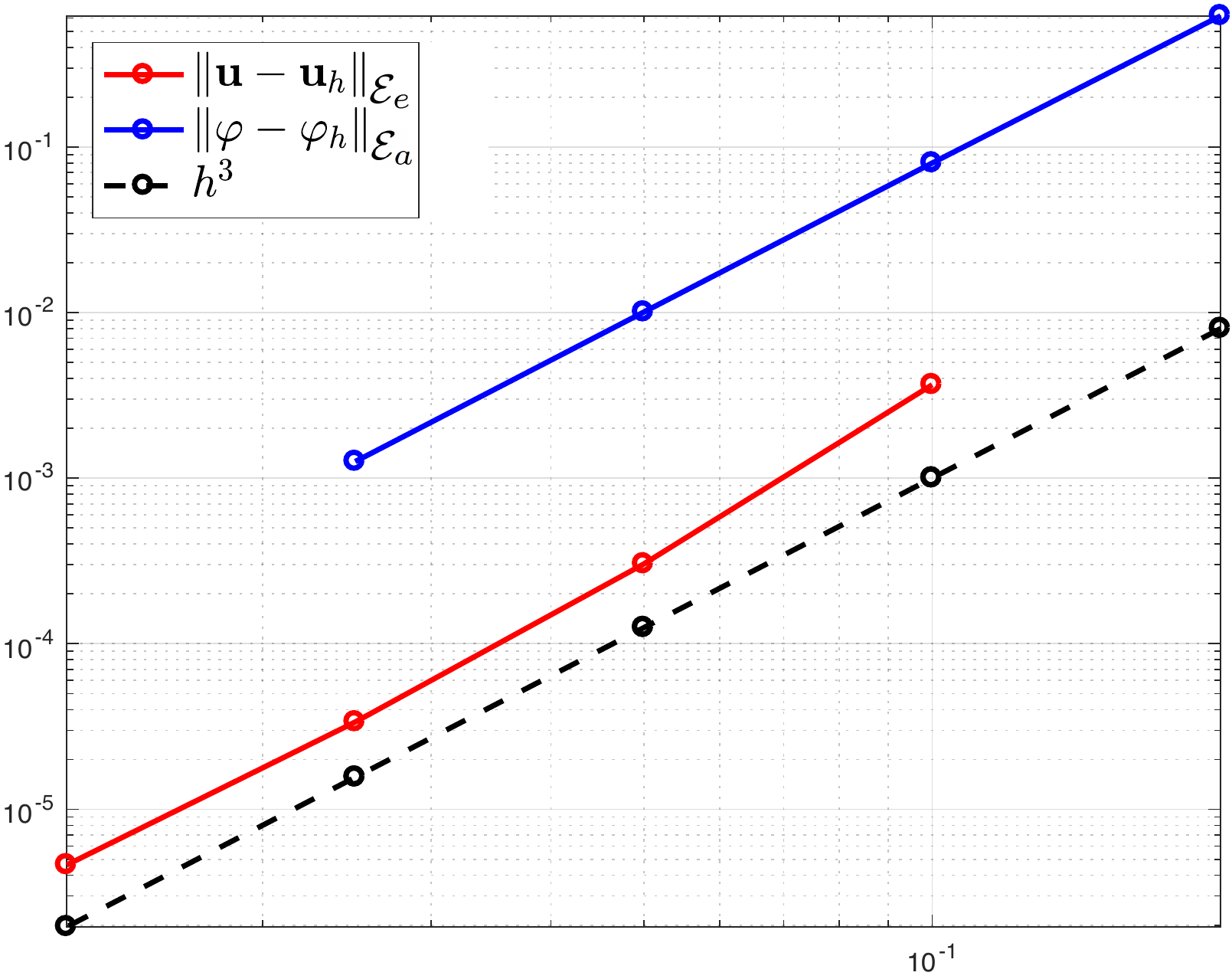}}\\
\subfloat[$N=N_e=N_a=2$]{\includegraphics[scale=.485]{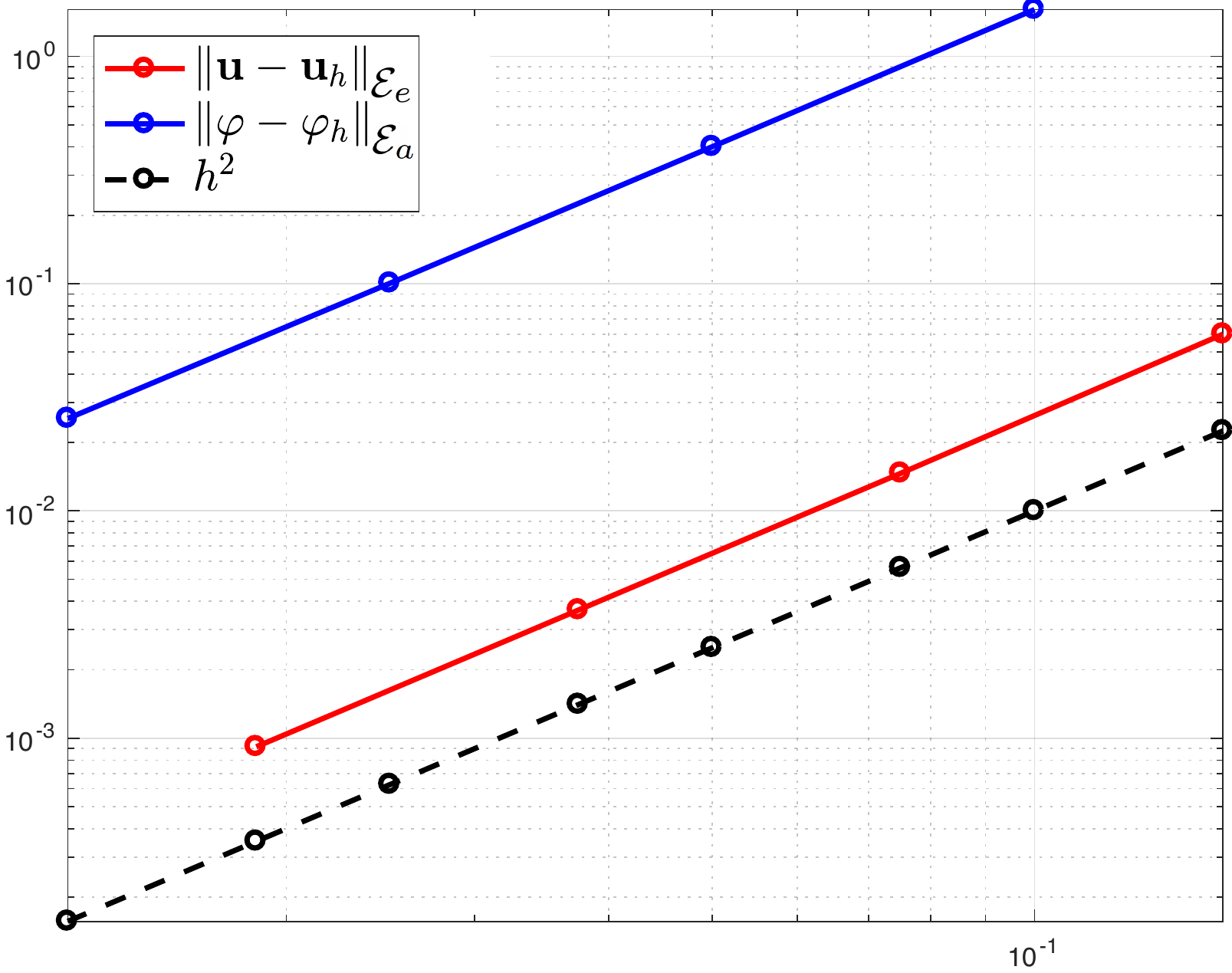}}
\subfloat[$h_e=0.1$, $h_a=0.2$, $N$ ranging from 2 to 7]{\includegraphics[scale=.48]{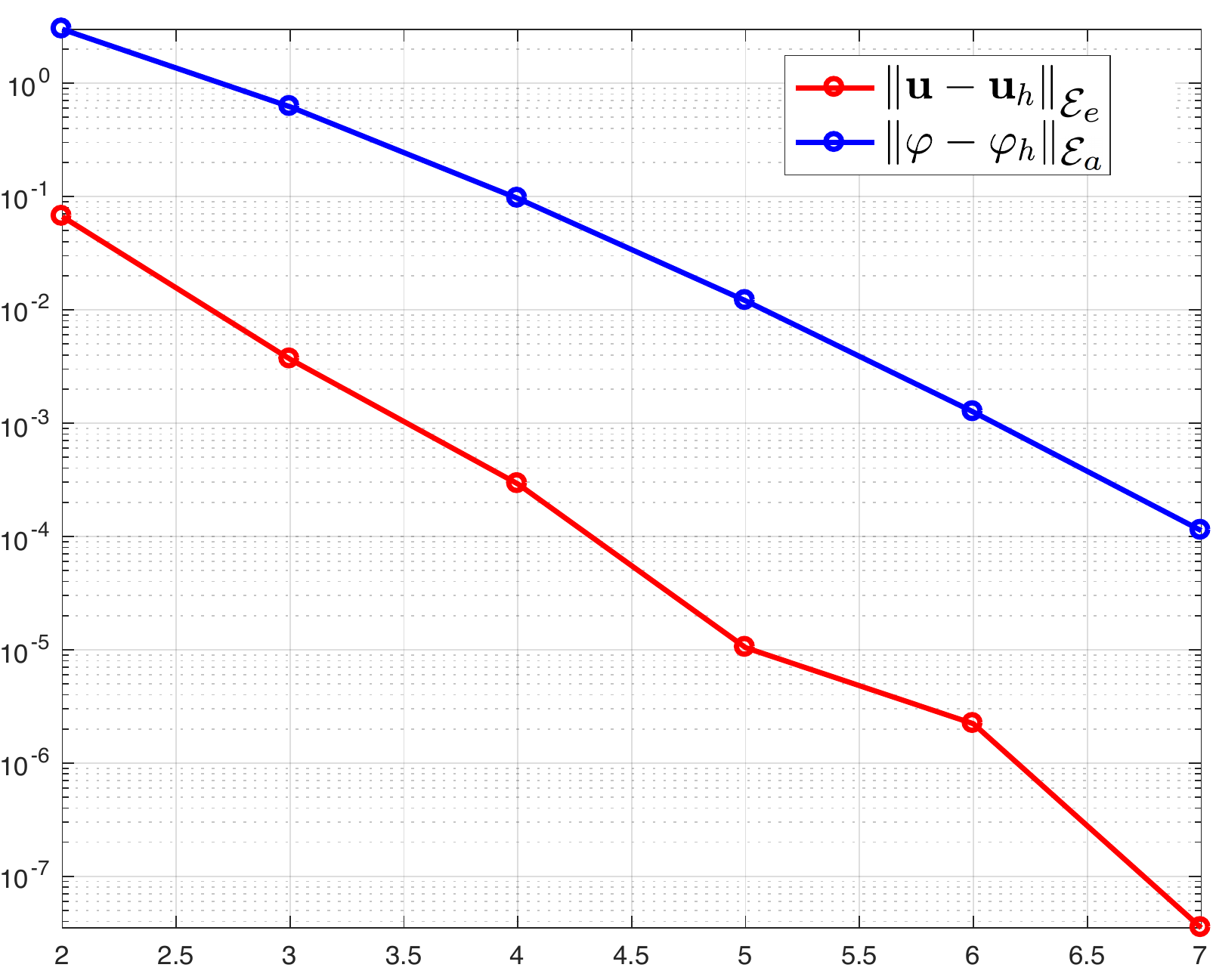}}
\caption{Test case~\ref{sec_test}. Error error vs.~$h$ (a)--(b)--(c) and $N$ (d) at $t=0.1\,s$ for non-matching hexahedral grids. Initial meshsizes are $h_e=0.1$, $h_a=0.2$ in (a) and (b), and $h_e=0.1$, $h_a=0.15$ in (c).}
\label{errore_h_NC}
\end{figure}

In this section we solve problem \eqref{strong_form}  in the parallelepiped $\Omega = (-1,1) \times (0,1) \times (0,1)$ on both matching and non-matching  grids (Figure~\ref{mesh}), and {verify the convergence results shown in Theorem~\ref{th-conve}}. Here $\Oe = (-1,0)\times (0,1) \times (0,1)$ and $\Oa = (0,1)^3$; the interface is thus given by $\Gamma_{\mr I} = \{0\}\times [0,1] \times [0,1]$. In all cases we compute the energy norm of the error at time $t = 0.1${, cf. \eqref{errore'}}. {For the time discretization we employed the staggered scheme presented in the previous section}. The timestep will be precised depending on the case under consideration. Finally, we choose $\rho_e = 2.7$, $c_P = 6.20$, $c_S = 3.12$, $\rho_a = 1$, and $c=1$ (cf.~\cite{monkola-sanna,monkola-sanna-paper}).
The right-hand sides $\bf f_e$ and $f_a$ are chosen so that the exact solution is given by
\begin{equation}\label{nonhom_bc}
\begin{aligned}
\bf u(x,y,z;t) & = \left(\cos\Big(\frac{4\pi x}{c_P}\Big), \, \cos\Big(\frac{4\pi x}{c_S}\Big), \, \cos\Big(\frac{4\pi x}{c_S}\Big)\right)\cos( 4\pi t), \\
\phi(x,y,z;t) & = \sin\Big(\frac{4\pi x}{c}\Big)\sin(4\pi t).
\end{aligned}
\end{equation}
Grids are sequentially refined starting from an initial mesh with uniform meshsize $h = 0.1$ in the matching case (Figure~\ref{mesh}a); on the other hand, in the non-matching case (Figure~\ref{mesh}b--\ref{mesh}c), the submeshes of $\Oe$ and $\Oa$ have the two initial respective meshsizes $h_e = 0.1$ and $h_a = 0.2$. Numerical tests carried out in both matching (Figure~\ref{errore_h}) and non-matching (Figure~\ref{errore_h_NC}) cases, show that $h$- and {$N$}-convergence rates match those predicted by {\eqref{errore'}}.
We also considered a further non-matching grid, where the initial submeshes are such that their meshsizes are not a multiple of each other, i.e., $h_e = 0.1$ and $h_a = 0.15$ (Figure~\ref{errore_h_NC}c). In this case, we obtain {a} quadratic {order of} convergence, {as expected with polynomial degree $N_e=N_a=2$}.
\subsection{Scholte waves}
\label{sec_scholte}
\begin{figure}
\includegraphics[keepaspectratio,scale=0.7]{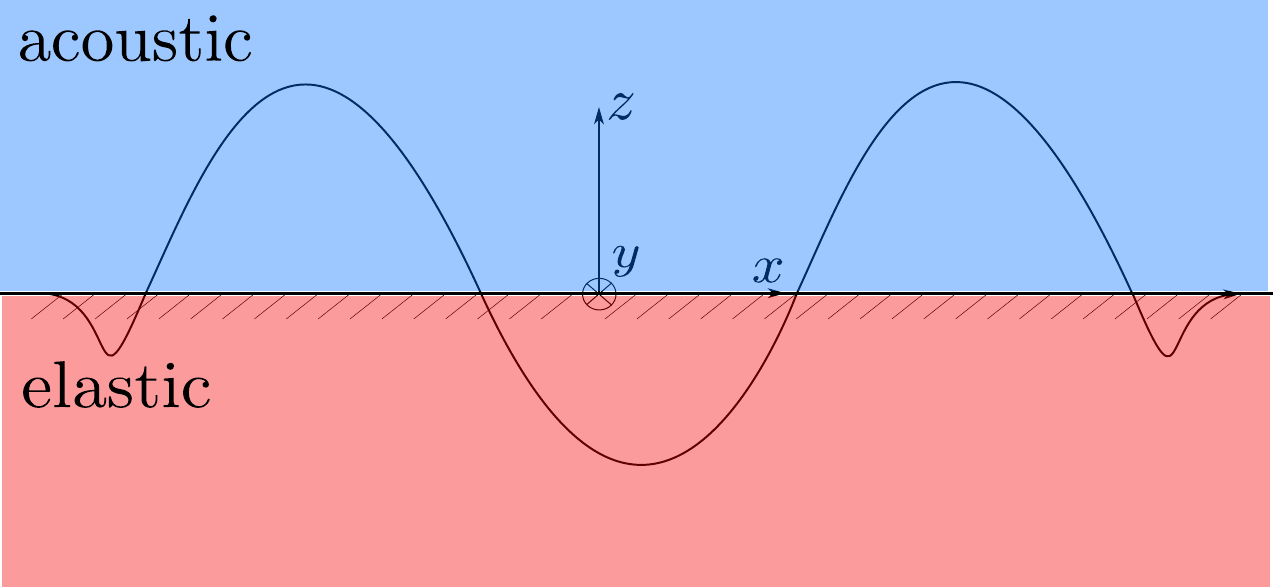}
\caption{Test case~\ref{sec_scholte}. Scholte wave at the interface between an elastic medium and an acoustic one.}
\label{scholte_fig}
\end{figure}

Scholte waves are an example of boundary waves, propagating along elasto-acoustic interfaces (cf.~Figure~\ref{scholte_fig}). Their amplitude decays exponentially away from the interface. As in \cite{ghattas}, we consider here two half-spaces. The lower half, $z<0$, is occupied by an elastic medium with $\lambda = \mu = 1$ and $\rho_e=1$; and the upper half, $z>0$, by an acoustic medium with $c=1$ and $\rho_a=1$. The analytic expressions of the displacement field $\bf u \equiv (u_1,u_2,u_3)$ and velocity potential $\phi$ can be inferred from \cite{ghattas}, where a displacement-based formulation is employed in both domains (see also \cite[Section~5.2]{kaufman-levshin}), and are the following. For $z < 0$ (elastic region), we have
\begin{equation}\label{scholte_e}
\begin{aligned}
 u_1(x,y,z;t) & =  k(B_2 e^{k  b_{2p} z} - 
  B_3 b_{2s} e^{k  b_{2s} z} ) \cos(k x - \omega t),\\
 u_2(x,y,z;t) & = 0, \\
  u_3(x,y,z;t) & =  k(B_2 b_{2p} e^{k b_{2p} z} - B_3 e^{k b_{2s} z })  \sin(k x - \omega t); 
  \end{aligned}
\end{equation}
and, for $z>0$ (acoustic region),
\begin{equation}\label{scholte_a}
\phi(x,y,z;t) = \omega B_1 e^{-k b_{1p} z} \cos(kx - \omega t).
\end{equation}
Here, the wavenumber is $k = \omega/\cs$, for a given frequency $\omega$ and Scholte wave speed $\cs$. The decay rates are given by
$$b_{1p} = \sqrt{1-\frac{\cs^2}{c^2}},\qquad b_{2p} = \sqrt{1-\frac{\cs^2}{c_P^2}}, \qquad b_{2s} = \sqrt{1-\frac{\cs^2}{c_S^2}}.$$
Wave amplitudes $B_1$, $B_2$, and $B_3$ have to satisfy a suitable eigenvalue problem, say \mbox{$\bm\Lambda \bf B = \bf 0$} with $\bm\Lambda$ a suitable $3\times3$ matrix and $\bf B \equiv [B_1\,\,  B_2\,\, B_3]^T$, stemming from the transmission conditions imposed on $\Gamma_{\mr I}$, i.e.~$\bm\sigma(\bf u)\bf n_e = -\rho_a\dot\phi\bf n_e$ and $\de\phi/\de\bf n_a = -\dot{\bf u}\dotp \bf n_a$. The value of the Scholte wave speed $\cs$ is thus given by the condition $\det \bm\Lambda = 0$. One can show that a Scholte wave speed exists for arbitrary combinations of material parameters. Based on the values of the material parameters we selected, we obtain, analogously to \cite{ghattas}, $\cs = 0.7110017230197$, and we choose $B_1 = 0.3594499773037$, $B_2 = 0.8194642725978$, and $B_3 = 1$. Also, for our numerical experiments, we choose $\omega = 1$, which gives, in turn, $k = 1.4064663525$.

\begin{figure}
\centering
\subfloat[]{\includegraphics[scale=.475]{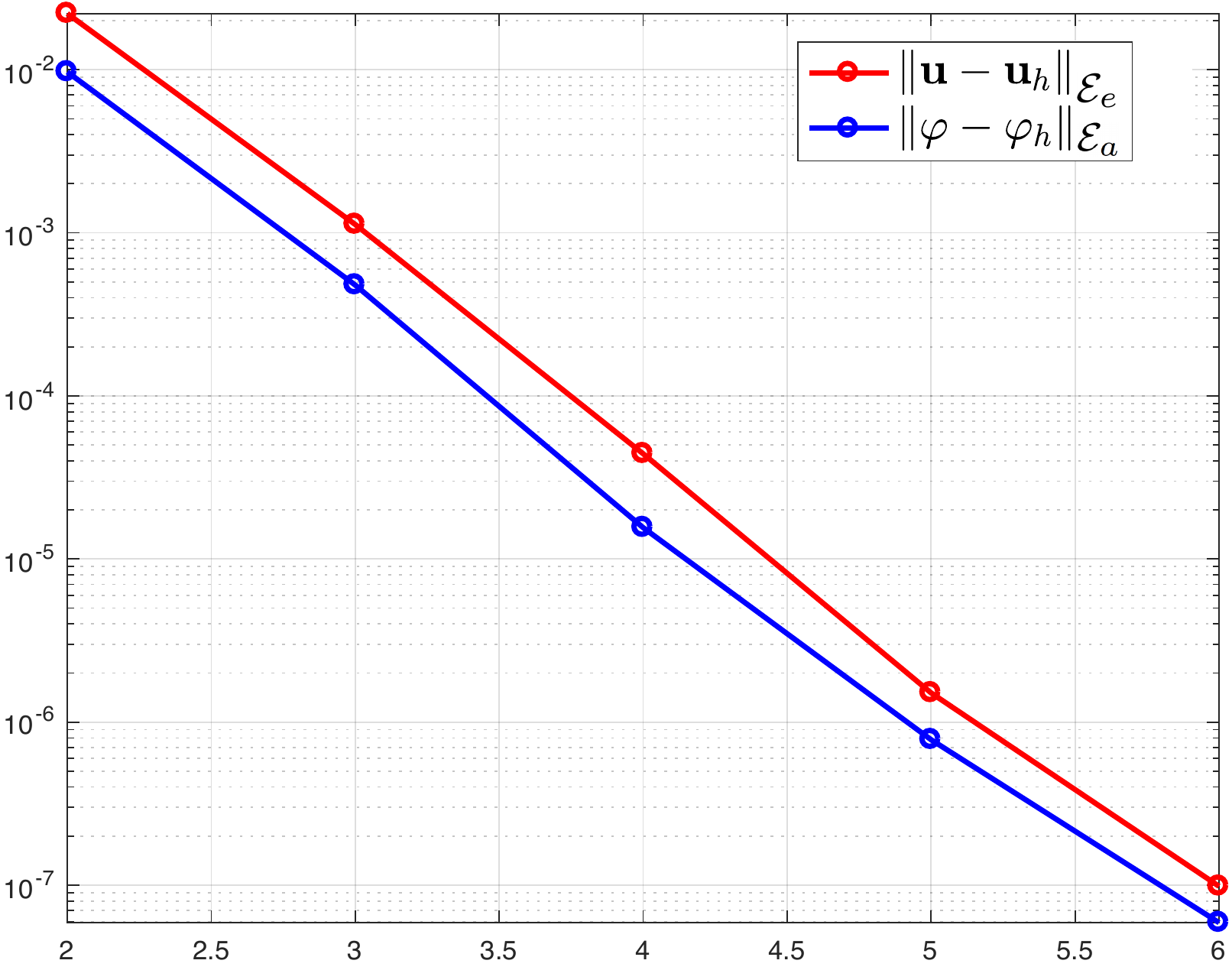}}
\subfloat[]{\includegraphics[scale=.475]{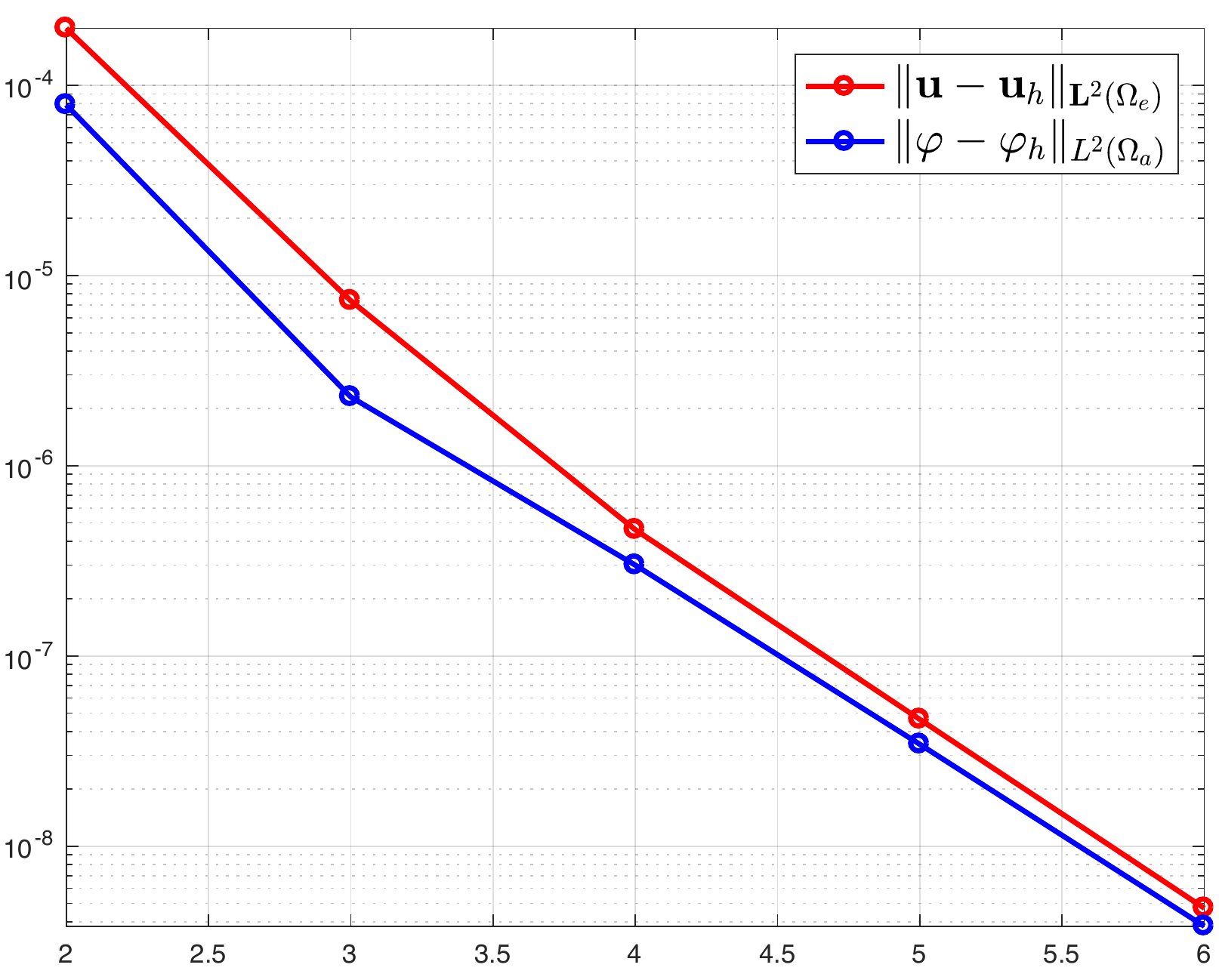}}
\caption{Test case~\ref{sec_scholte}. Error in the energy (a) and $L^2$ (b) norms vs.~$N$ at $t=0.1\,s$, with $N$ ranging from 2 to 6.}
\label{scholte_conv}
\end{figure}

We use a uniform mesh consisting of 2400 elements (corresponding to a meshsize $h = h_e = h_a = 0.416$) over the domain $(-1, 1) \times ( -1, 1) \times  (- 20, 20)$, and we impose Dirichlet conditions all over the boundary. Figure \ref{scholte_conv} shows asymptotic exponential convergence rate of the error in the {energy} and $L^2$ norms, as expected.
%
\subsection{Underground acoustic cavity}
\label{sec_cavity}
\begin{figure}
\centering
\includegraphics[scale=.175]{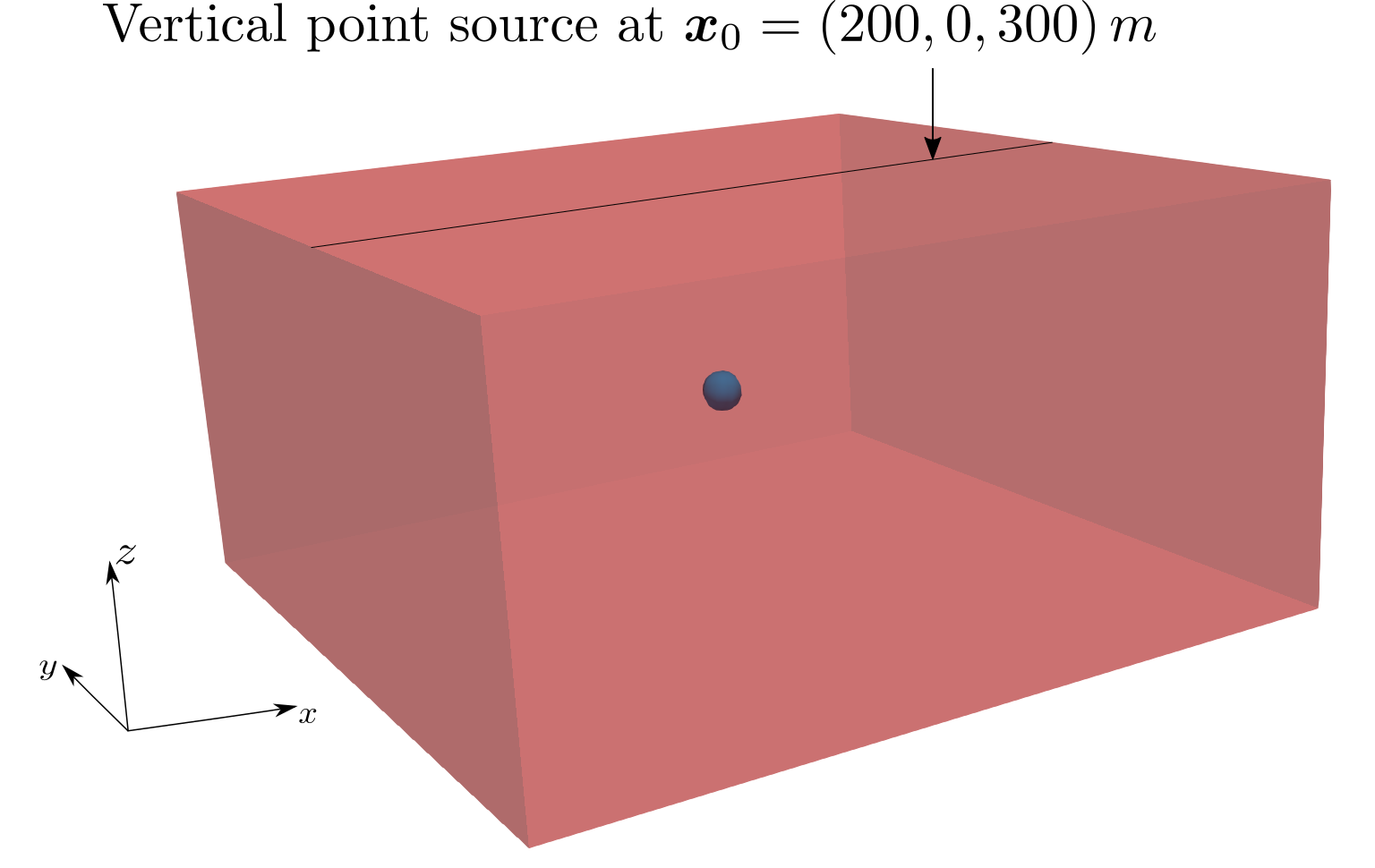}
\caption{Test case~\ref{sec_cavity}. Geometry of the computational domain for the case of a seismic wave in the presence of an underground cavity.}
\label{cavita}
\end{figure}

As a last test case, we simulate a seismic wave in the presence of an underground spherical acoustic cavity. This problem arises in several applications, the most important one, besides non-destructive testing \cite{phd-ferroni}, is given by near-surface seismic studies to detect the presence of cavities in the subsoil, which are originated after underground nuclear explosions, and can give rise to resonance effects when a seismic event occurs \cite{perugia-cavity-2}. In particular, the geometry we consider is the following: the acoustic domain is given by an open ball $\Omega_a = \{\bm x\in \bb R^3 : \|\bm x\| < R\}$, of radius $R=30\,m$, and the elastic one is $\Omega_e = (-L_x,L_x)\times(-L_y,L_y)\times(-L_z,L_z)\setminus\overline{\Omega}_a$ surrounding the cavity, with $L_x=L_y=600\,m$ and $L_z=300\,m$ (Figure~\ref{cavita}). Non-reflecting boundary conditions are imposed on the external elastic boundary. The system is excited by a point Ricker wavelet of the following form: $$\bf f_e(\bm x,t) = f(t)\bf e_z \delta(\bm x - \bm x_0),\quad f(t) = f_0 \left(1-2\pi^2 f_p^2(t-t_0)^2\right)e^{-\pi^2 f_p^2(t-t_0)^2},$$
with $\bf e_z = (0,0,1)$, $\bm x_0 = (200, 0, 300)\,m$, $t_0=0.25\,s$, $f_0=10^{10}\,N$, and \emph{peak frequency} $f_p$. The set of data and space discretization parameters is summarized in Table \ref{table_cavity}, where we write $c_P$ for $c$ in the case of an acoustic wave.
\begin{table}
 \caption{Test case~\ref{sec_cavity}. Material properties.}
  \label{table_cavity}{ 
\begin{tabular}{|c|c|c|c|}
  \hline
  Region  &  $\rho_{e/a}\, (kg/m^3)$ & $c_P\,(m/s)$ & $c_S\,(m/s)$ \\ \hline
  ${\Omega_e}$  & ${2700}$ & ${3000}$ & ${1734}$  \\ \hline
  ${\Omega_a}$ & ${1024}$ & ${300}$ & {--} \\ \hline
  \end{tabular}
  }
\end{table}
Since the wavelength inside the cavity is much smaller than outside, we are led to choosing a finer meshsize inside the cavity, and thus employ the following meshsizes: $h_e = 20\,m$, $h_a = 5\,m$. 
We use a polynomial degree $N_e=N_a=4$ on both domains, and we set the time-step to $\Delta t = 10^{-5}\,s$. 

\begin{figure}
\centering
\subfloat[]{\includegraphics[scale=.155]{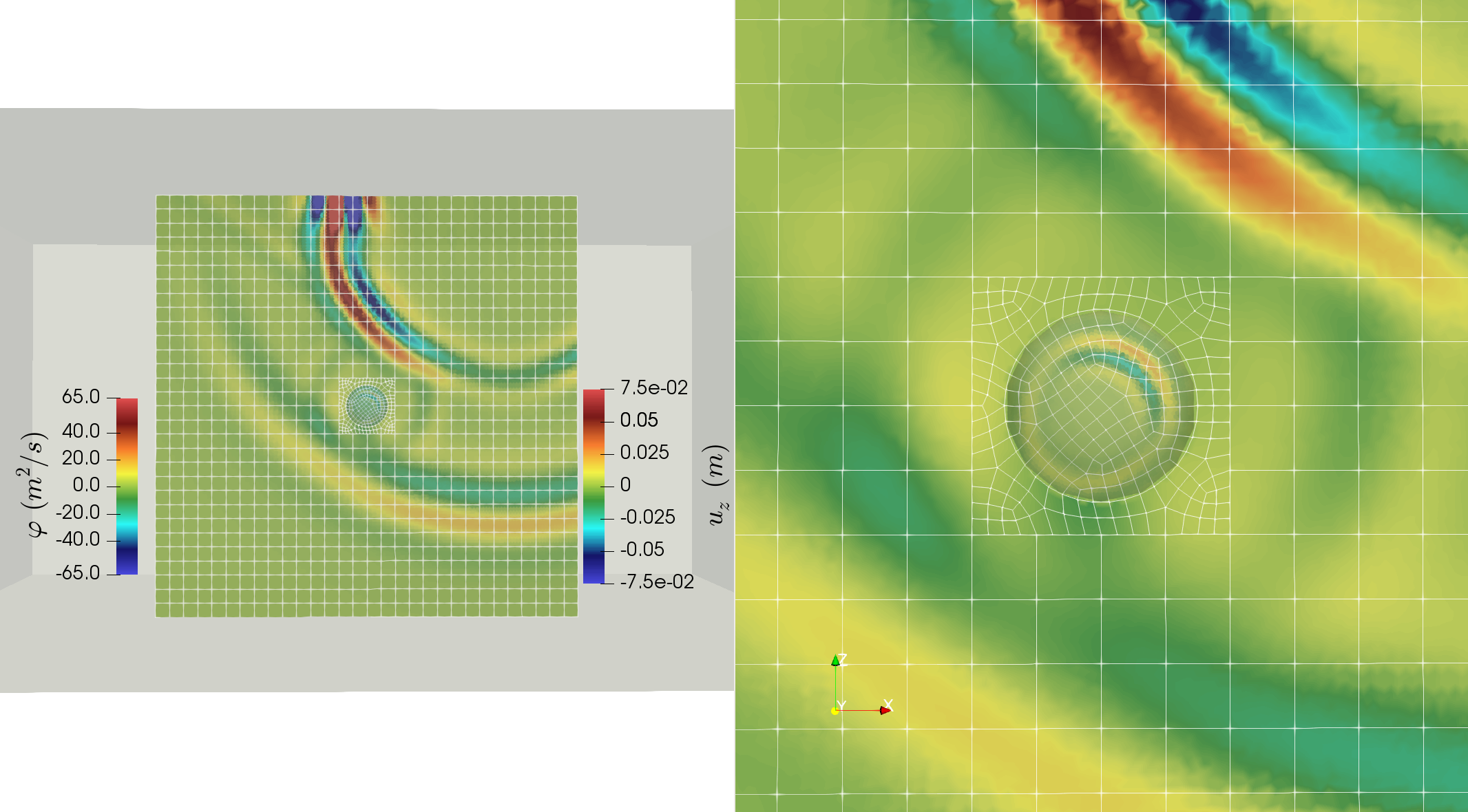}}
\hspace{.01cm}
\subfloat[]{\includegraphics[scale=.155]{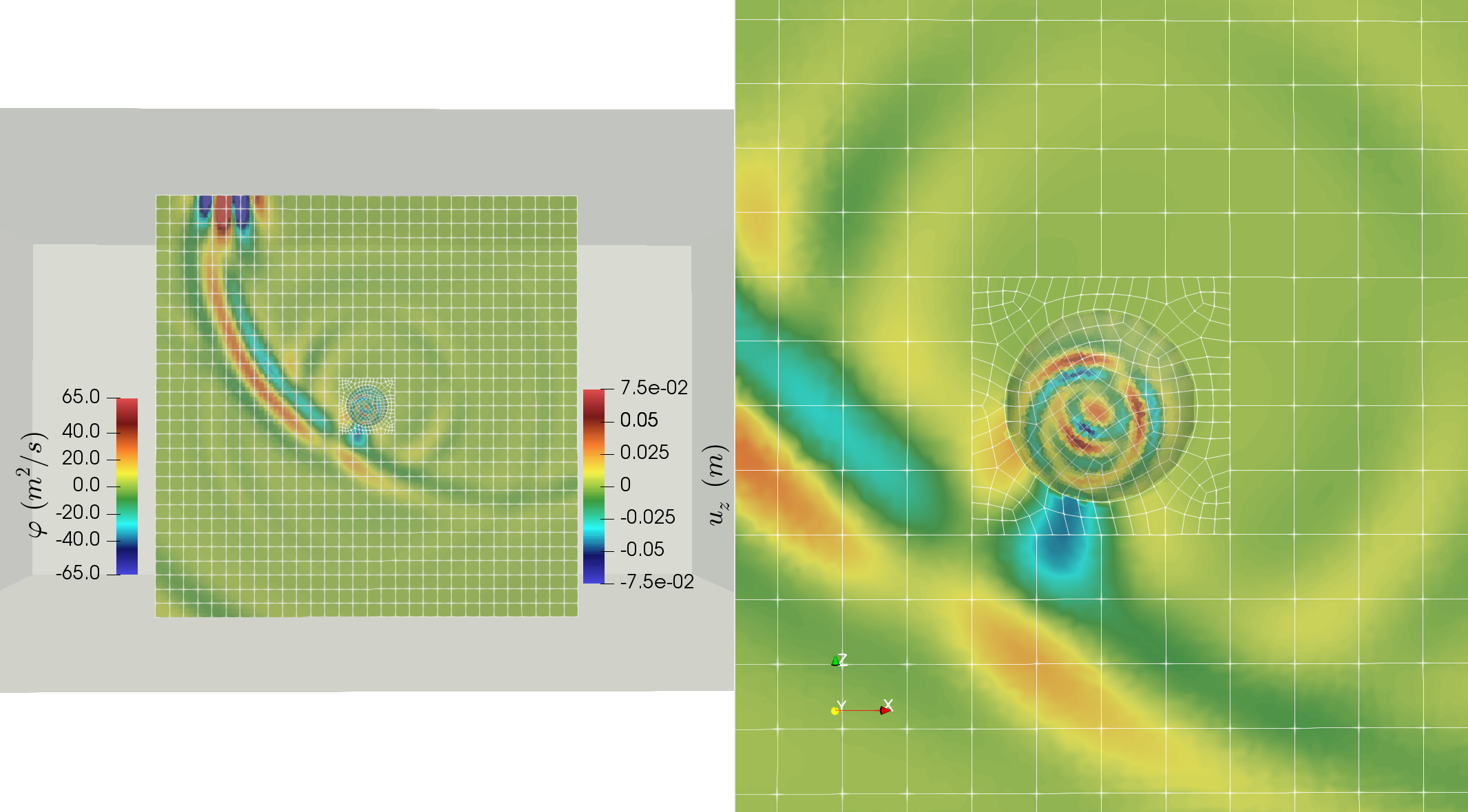}}
\hspace{.01cm}
\subfloat[]{\includegraphics[scale=.155]{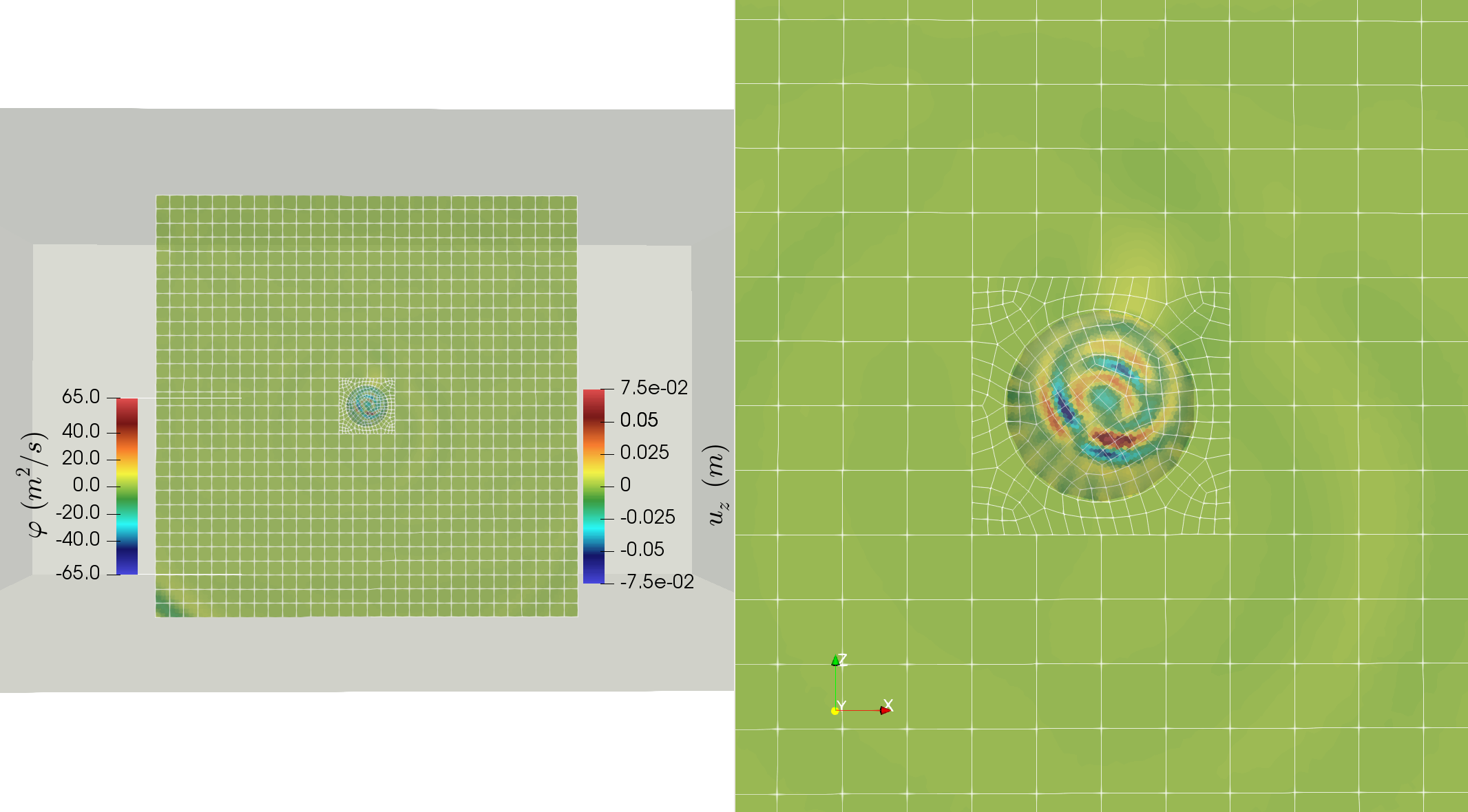}}
\caption{Test case~\ref{sec_cavity}. Displacement along the $z$-direction and velocity potential at time $t=0.4\,s$ (a), $t=0.5\,s$ (b), and $t=0.7\,s$ (c), for $f_p=22\,Hz$.}
\label{cavita1}
\end{figure}

\begin{figure}
\centering
\subfloat[]{\includegraphics[scale=.155]{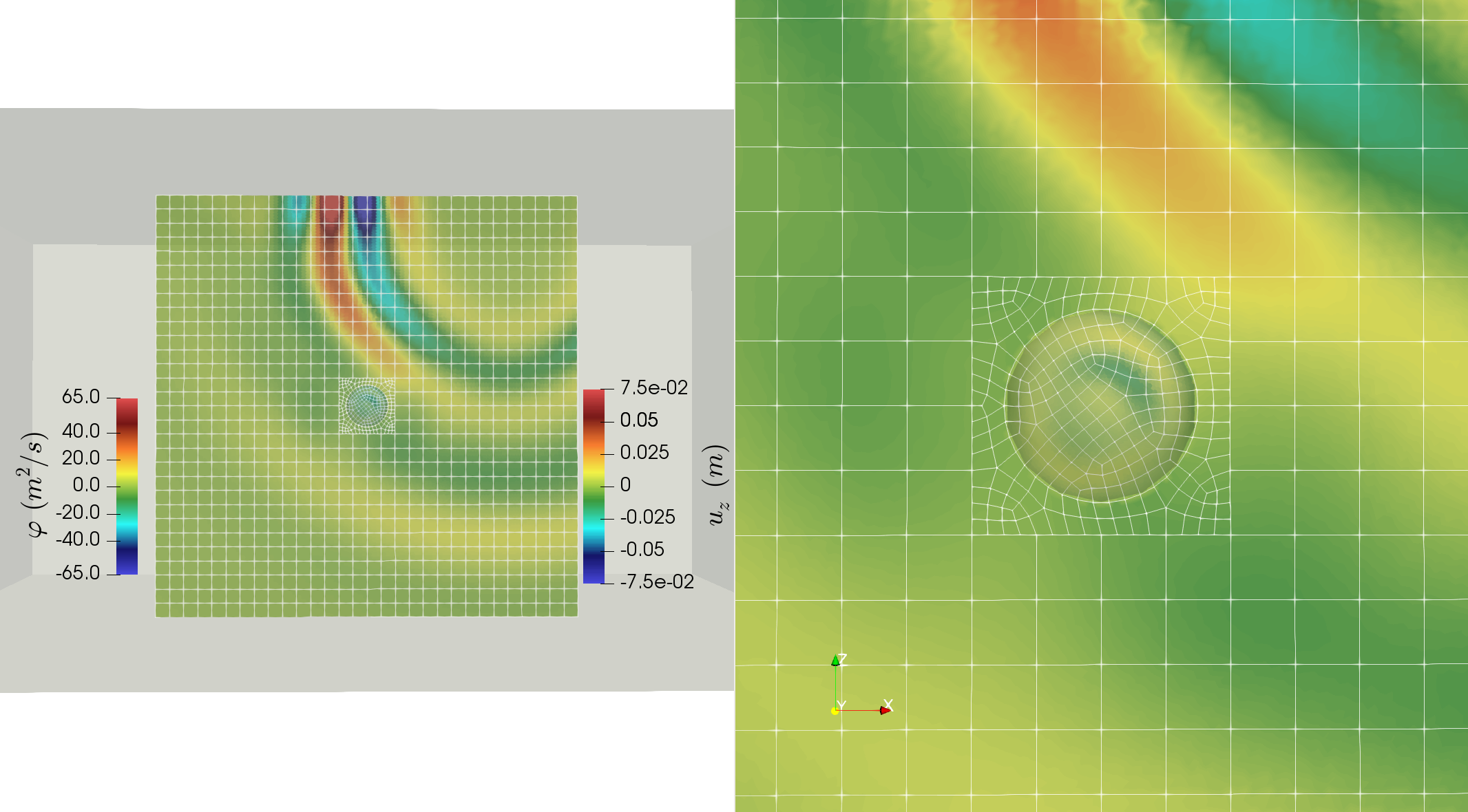}}
\hspace{.01cm}
\subfloat[]{\includegraphics[scale=.155]{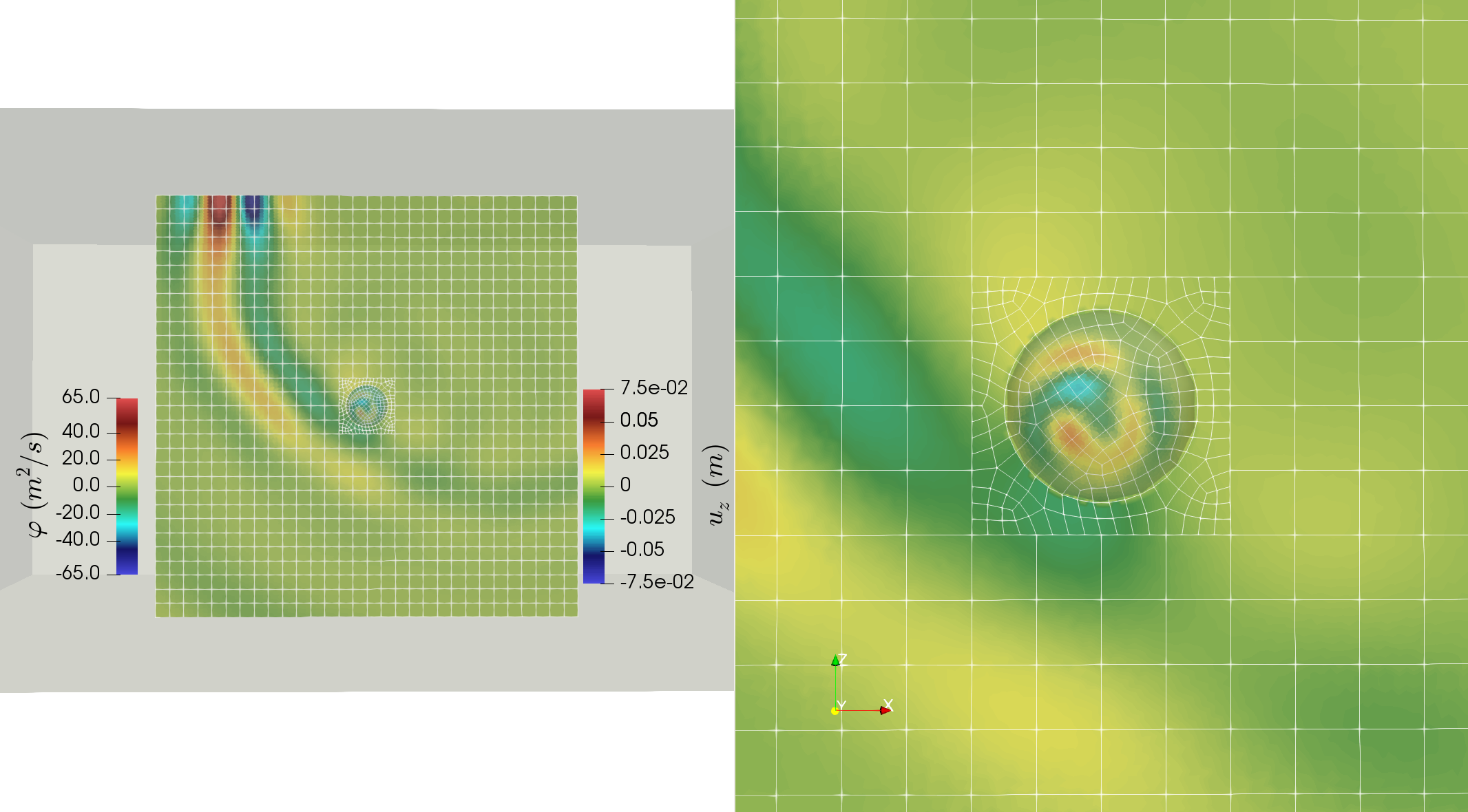}}
\hspace{.01cm}
\subfloat[]{\includegraphics[scale=.155]{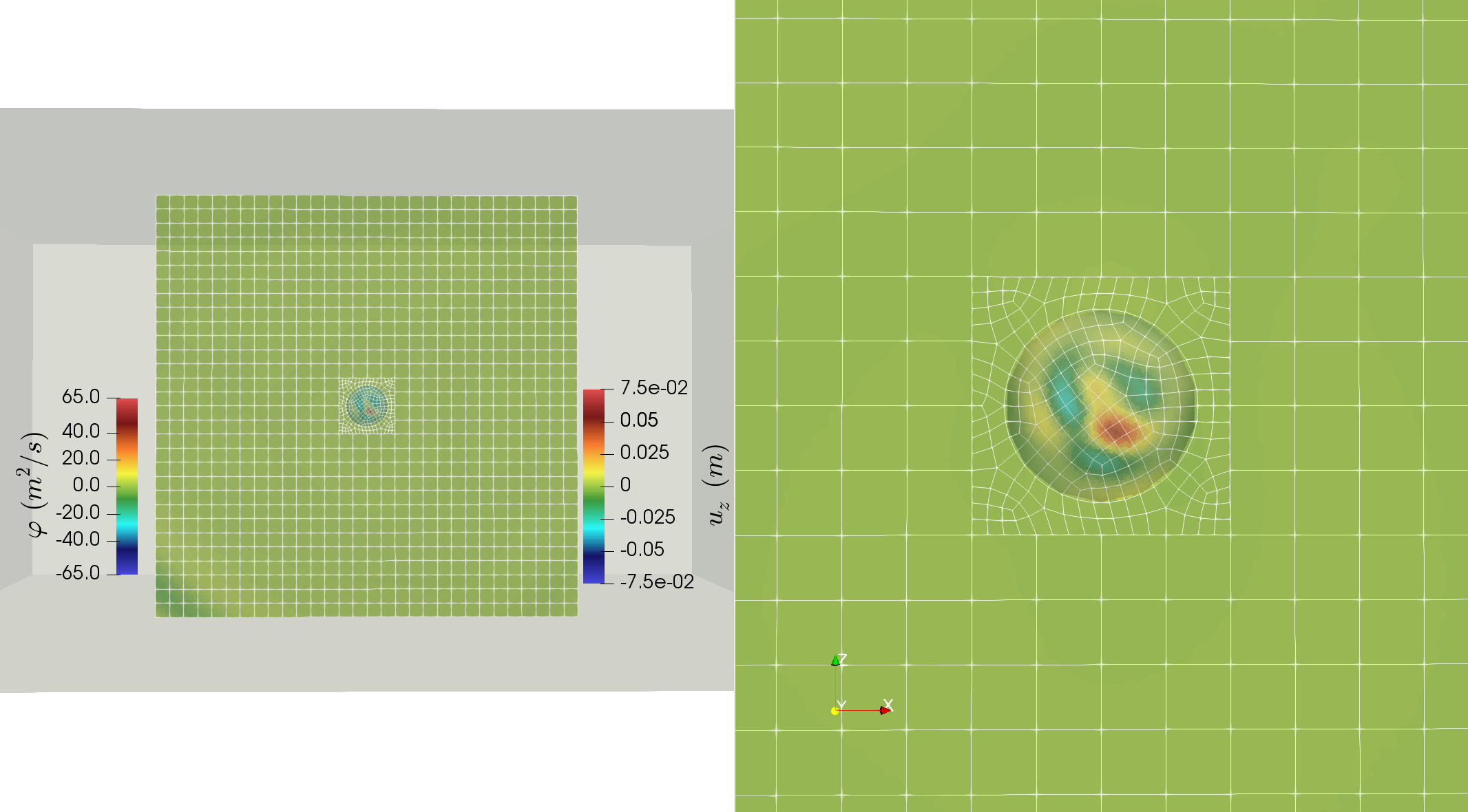}}
\caption{Test case~\ref{sec_cavity}. Displacement along the $z$-direction and velocity potential at time $t=0.4\,s$ (a), $t=0.5\,s$ (b), and $t=0.7\,s$ (c), for $f_p=11\,Hz$.}
\label{cavita2}
\end{figure}

\begin{figure}
\includegraphics[keepaspectratio,scale=0.6]{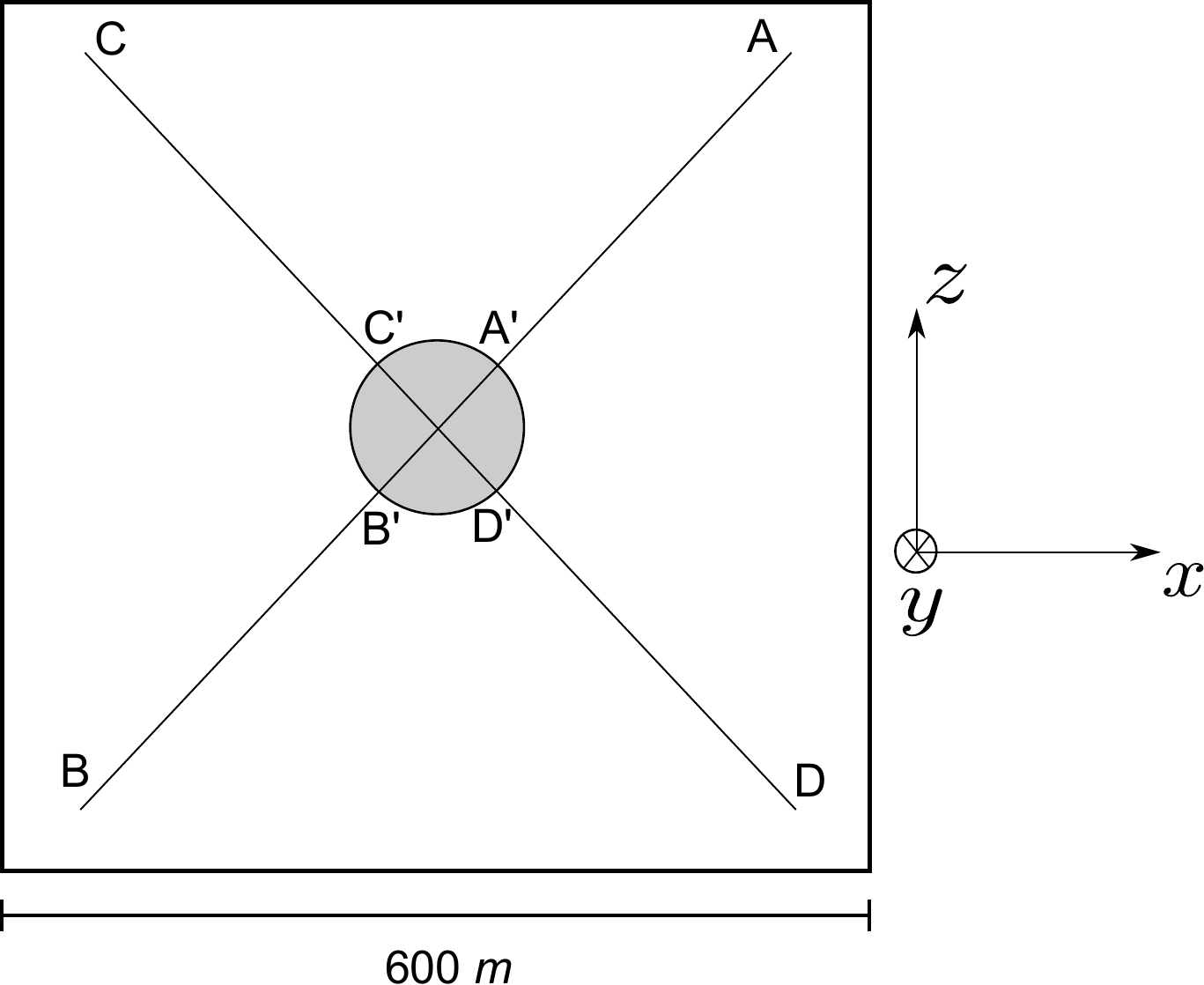}
\caption{Test case~\ref{sec_cavity}. Set of monitors in the square cross section of the computational domain {lying in the $xz$-plane, centered in the origin, with side $600\,m$}.}
\label{monitors}
\end{figure}

\begin{figure}
\centering
\subfloat{\includegraphics[scale=.2]{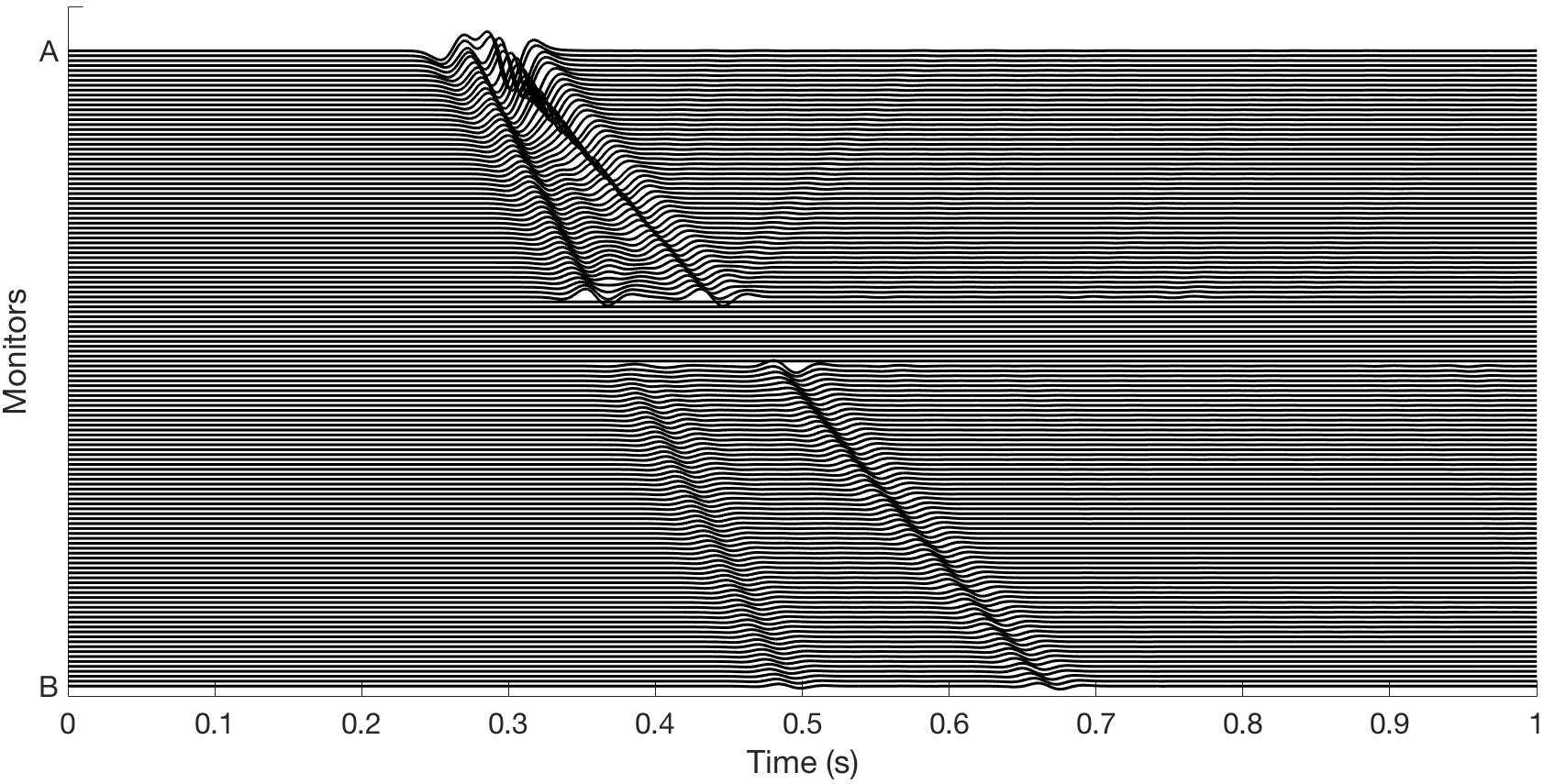}}
\hspace{.01cm}
\subfloat{\includegraphics[scale=.2]{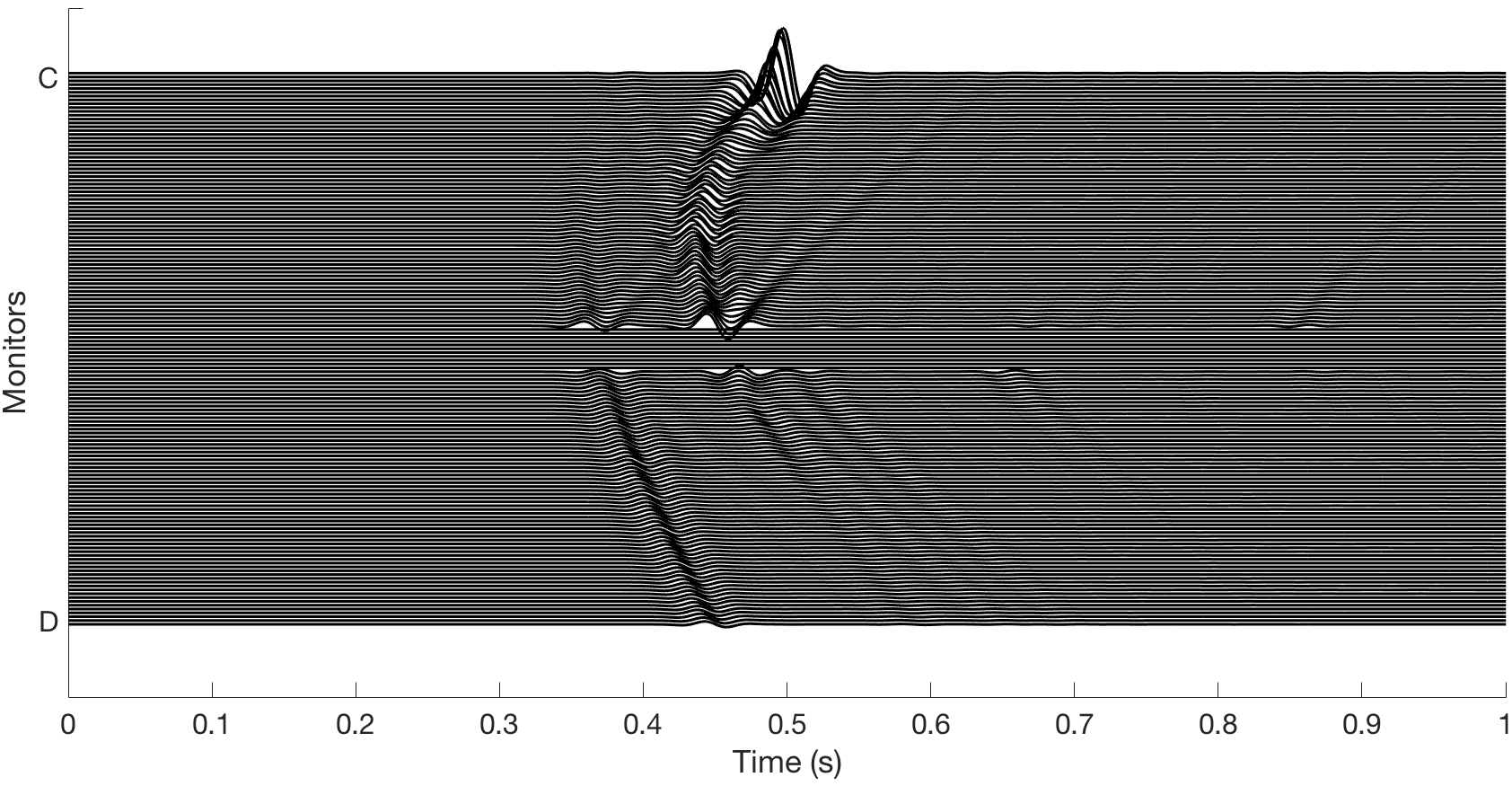}}
\caption{Test case~\ref{sec_cavity}. Time histories of the displacement along the $z$-direction for the monitored points in the elastic subsoil, for $f_p=22\,Hz$.}
\label{punti_x_1_e}
\end{figure}

\begin{figure}
\centering
\subfloat{\includegraphics[scale=.2]{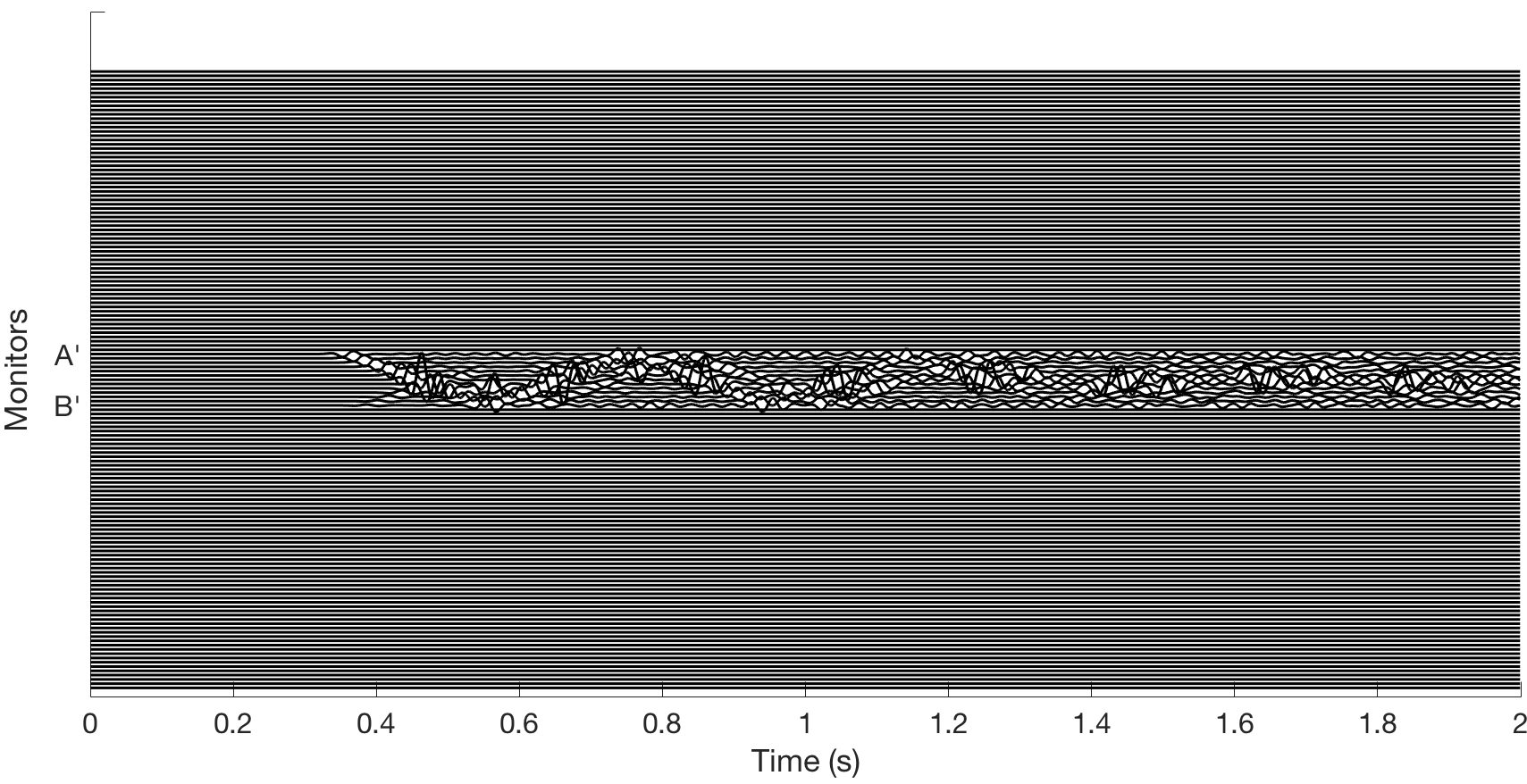}}
\hspace{.01cm}
\subfloat{\includegraphics[scale=.2]{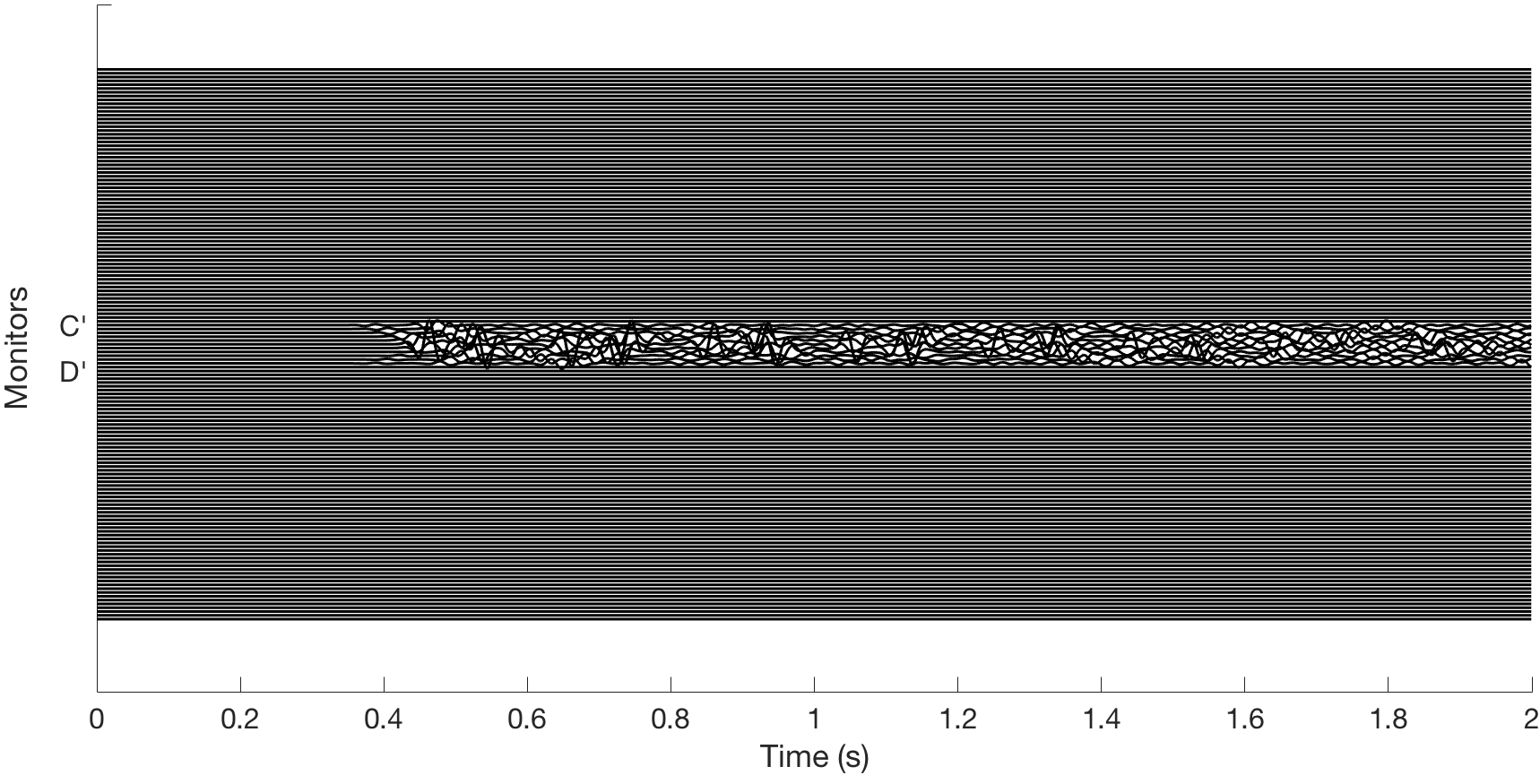}}
\caption{Test case~\ref{sec_cavity}. Time histories of the velocity potential for the monitored points in the acoustic cavity, for $f_p=22\,Hz$.}
\label{punti_x_1_ac}
\end{figure}

\begin{figure}
\centering
\subfloat{\includegraphics[scale=.2]{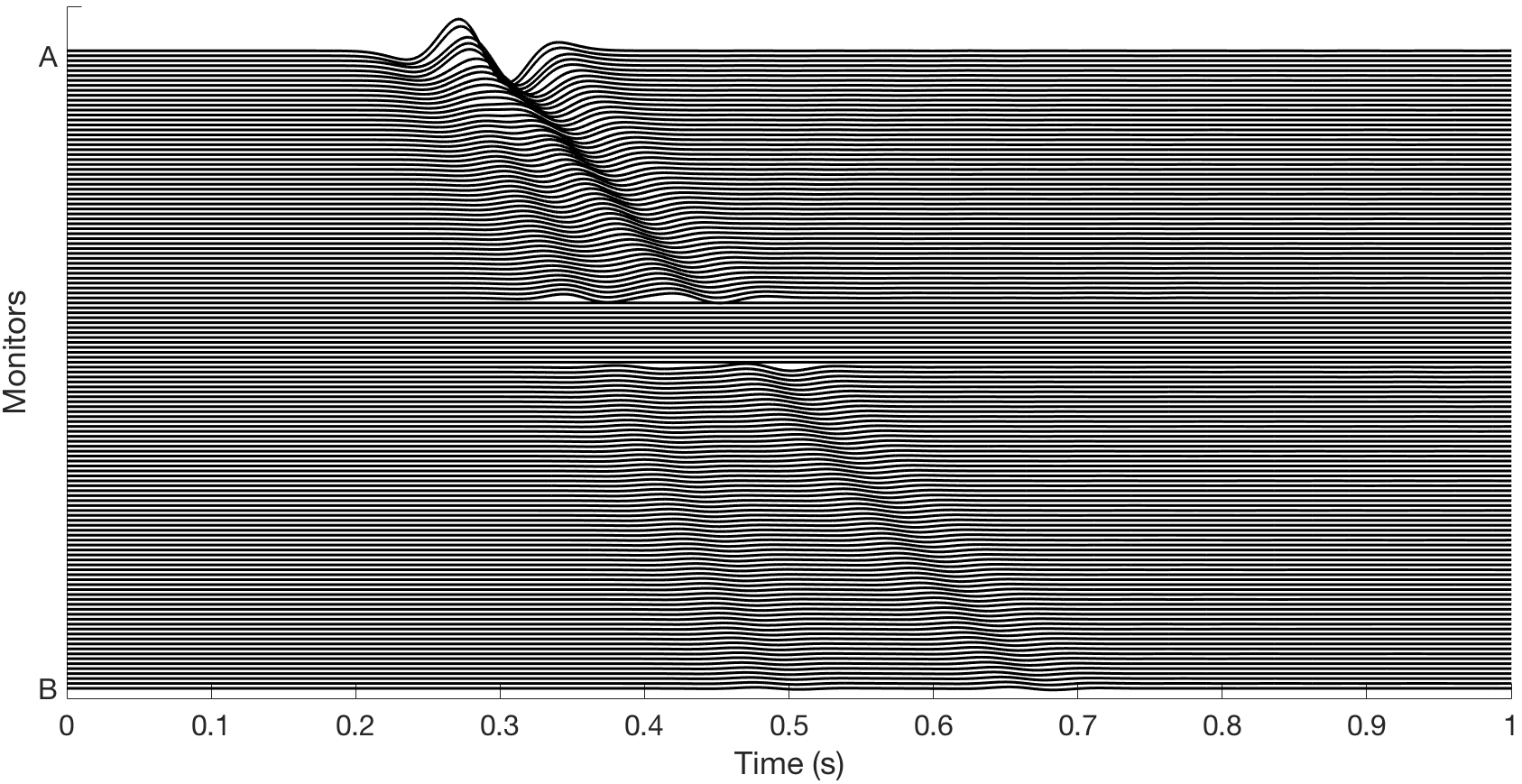}}
\hspace{.01cm}
\subfloat{\includegraphics[scale=.2]{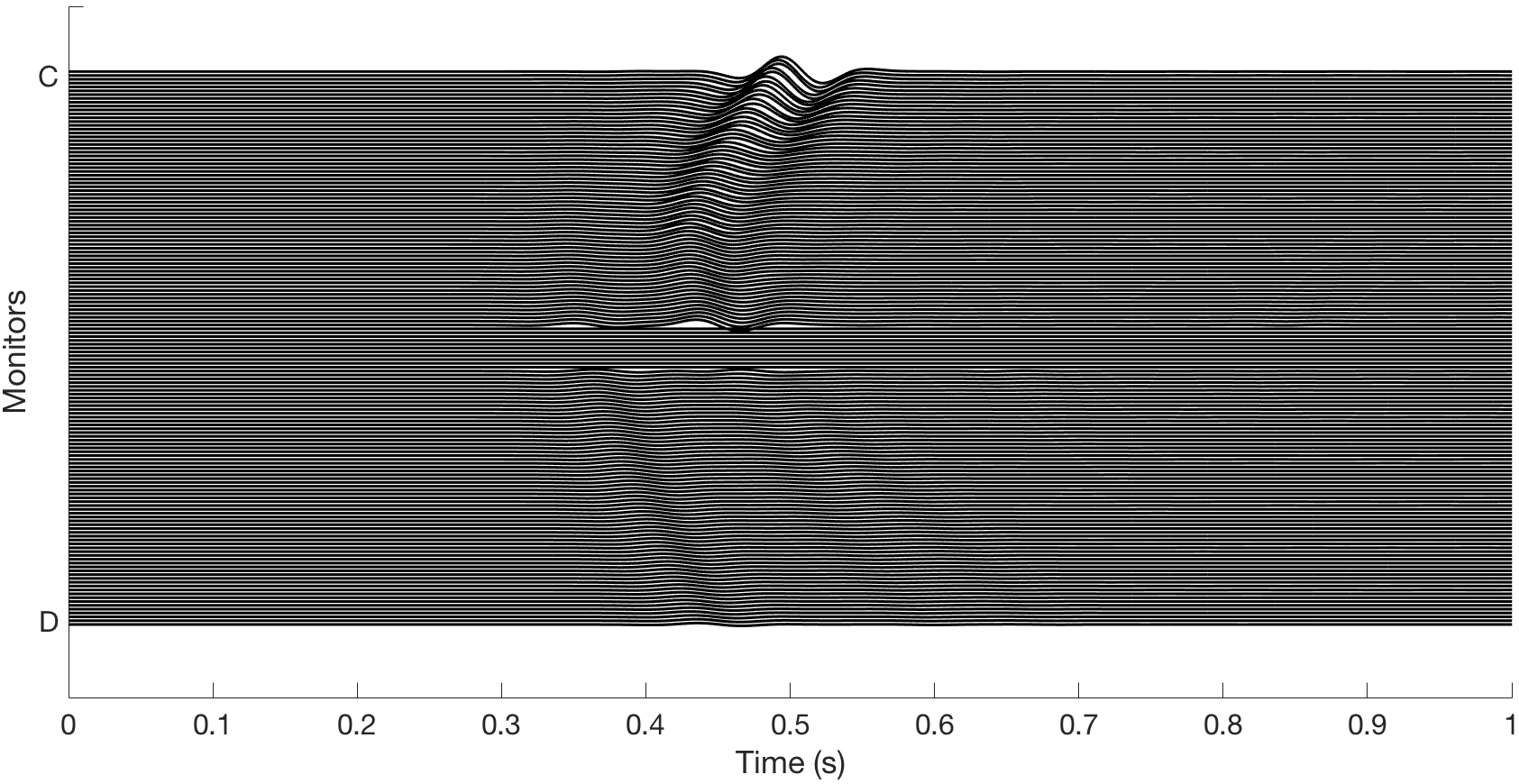}}
\caption{Test case~\ref{sec_cavity}. Time histories of the displacement along the $z$-direction for the monitored points in the elastic subsoil, for $f_p=11\,Hz$.}
\label{punti_x_2_e}
\end{figure}

\begin{figure}
\centering
\subfloat{\includegraphics[scale=.2]{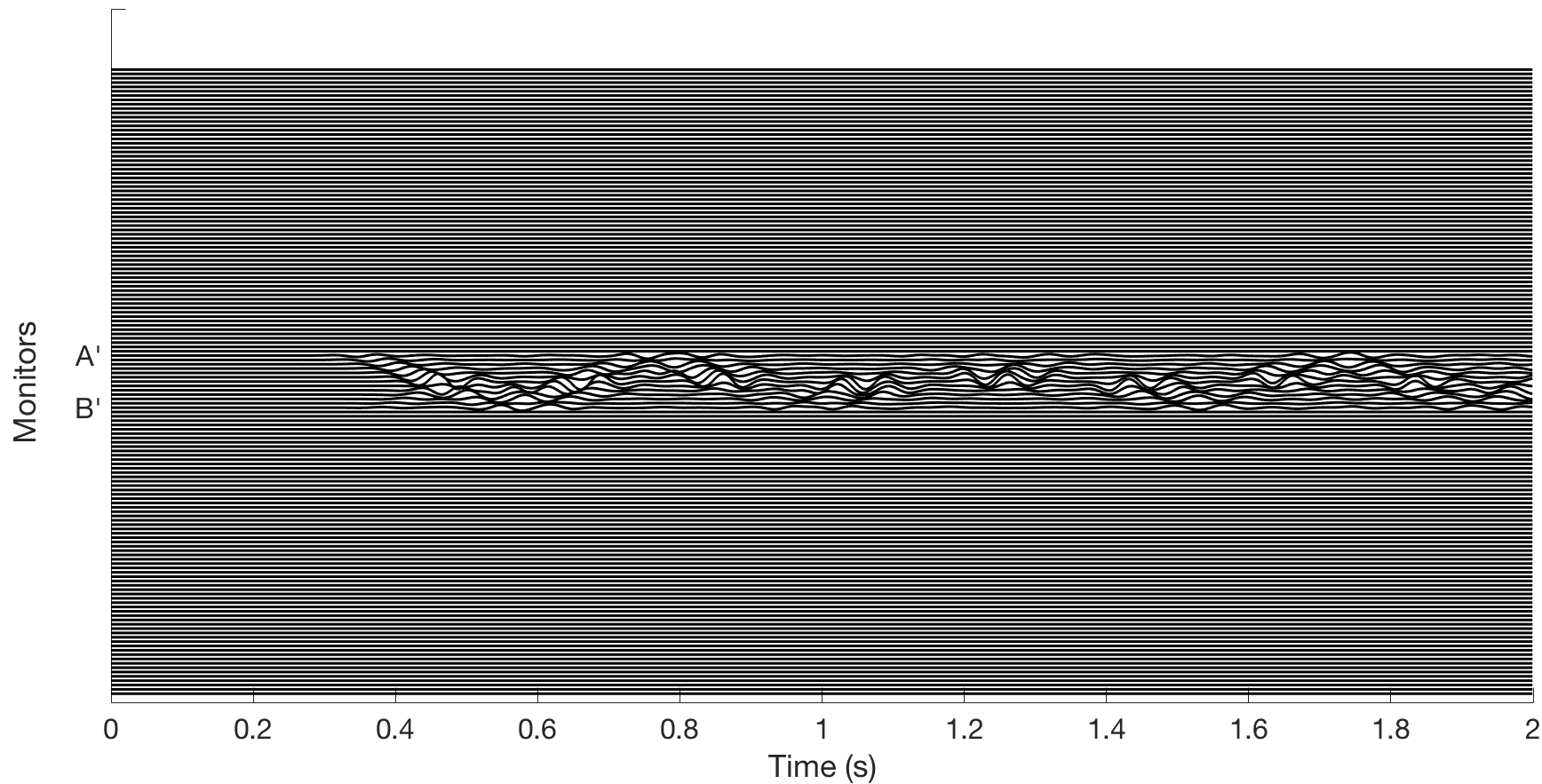}}
\hspace{.01cm}
\subfloat{\includegraphics[scale=.2]{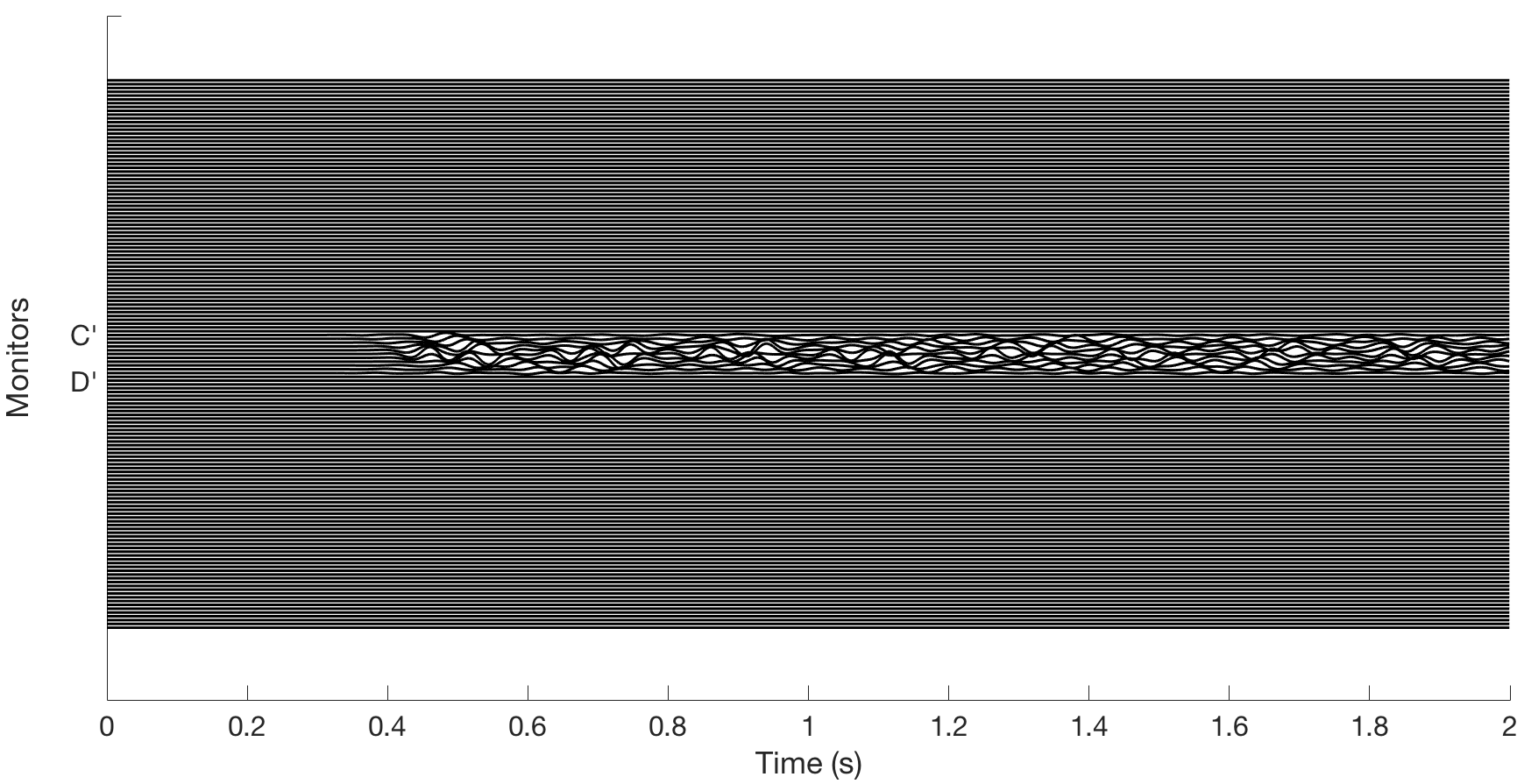}}
\caption{Test case~\ref{sec_cavity}. Time histories of the velocity potential for the monitored points in the acoustic cavity, for $f_p=11\,Hz$.}
\label{punti_x_2_ac}
\end{figure}

Figure~\ref{cavita1} shows the $z$-component $u_z$ of the displacement field in the subsoil and the acoustic velocity potential $\phi$ in the spherical cavity at times $t=0.4\,s$, $t=0.5\,s$, and $t=0.7\,s$ when the peak frequency is set to $f_p=22\,Hz$, whereas Figure~\ref{cavita2} shows the same quantities when $f_p=11\,Hz$.
We remark that, in the first case (Figure~\ref{cavita1}), the elastic wave detects the acoustic cavity: spherical wavefronts are generated due to refraction phenomena between the cavity and the subsoil, since the wavelength corresponding to the value $f_p = 22\,Hz$ is comparable with the diameter of the cavity. On the other hand, if the peak frequency is reduced by a factor two (Figure~\ref{cavita2}), we observe that the interaction of the elastic wave with the cavity is weaker than in the first case, since the corresponding wavelength is twice as much as in the first case. In both cases, since outside the sphere the material is stiff, the acoustic wave remains trapped within the cavity over time and it generates reflection and refraction effects.
These phenomena can be better represented and remarked if the time histories of a number of monitored points in the elastic and acoustic domains are considered. In particular, we took into account an $X$-shaped set of points in {a square cross section of the computational domain lying in the $xz$-plane, centered in the origin, with side $600\,m$} (Figure~\ref{monitors}). Time histories of points in the subsoil and in the underground cavity are showcased in Figures~\ref{punti_x_1_e} and \ref{punti_x_1_ac} for the first case ($f_p=22\,Hz$) and in Figures~\ref{punti_x_2_e} and \ref{punti_x_2_ac} for the second case ($f_p=11\,Hz$). In particular, reflection phenomena for elastic waves are clearly more remarkable in the first case than in the second. As expected, point A being the closest one to the location of the seismic source, is the first to undergo a displacement impulse, which is then delayed for the other points; the same occurs in the second case. Finally, in both cases, we clearly see that the acoustic wave remains trapped in cavity over time, due to persistent reflections.
\section{Conclusions and perspectives}
We have presented a Discontinuous Galerkin Spectral Element method for the approximation of the elasto-acoustic evolution problem. Several numerical experiments carried out in a three-dimensional framework have been discussed, both to {verify} the theoretical results and to simulate a scenario of physical interest. {Our approach is well-suited to comply with the requirements for the discretization of heterogeneous seismic wave propagation problems (geometric flexibility, high-order accuracy, and flexibility); in addition, it allows for the treatment of non-matching grids at the interface between the elastic and the acoustic domains, which can therefore be generated independently on each of the domains.}
All numerical experiments have been carried out using {the computer code} \texttt{SPEED} \cite{speed}, {freely available at {\tt http://speed.mox.polimi.it}}.

A future work consists in the extension to {general} polyhedral meshes in \texttt{SPEED}, in order to tame the computational cost of mesh generation and {enhance the} geometrical flexibility {of the numerical  discretization}. {As a} second perspective is given by the enrichment of the models.
%
%
%
%
%
%

%
\FloatBarrier

\end{document}